\documentclass[letterpaper,11pt]{article}

\usepackage{etoolbox}
\usepackage{xspace}
\usepackage[usenames,dvipsnames,svgnames,table]{xcolor}
\usepackage[colorlinks=true,linkcolor=blue,citecolor=ForestGreen]{hyperref}
\usepackage{amsmath,amsthm,amssymb}
\usepackage{enumerate}
\usepackage{pgfplots}
\usepackage{graphicx}
\usepackage{scrextend}
\usepackage[utf8]{inputenc}
\usepackage[T1]{fontenc}

\usepackage{paralist}

\usepackage{amsmath}
\usepackage{amssymb}
\usepackage{bm}
\usepackage{thmtools}
\usepackage{hyperref}
\usepackage[capitalise]{cleveref}
\usepackage{float}
\usepackage[color=green!40,textsize=small]{todonotes}
\usepackage{xspace}
\usepackage{orcidlink}
\usepackage[margin=1in]{geometry} 

\usepackage{setspace}
\allowdisplaybreaks

\Crefname{figure}{Figure}{Figures}
\Crefname{claim}{Claim}{Claims}

\crefformat{equation}{\textup{#2(#1)#3}}
\crefrangeformat{equation}{\textup{#3(#1)#4--#5(#2)#6}}
\crefmultiformat{equation}{\textup{#2(#1)#3}}{ and \textup{#2(#1)#3}}
{, \textup{#2(#1)#3}}{, and \textup{#2(#1)#3}}
\crefrangemultiformat{equation}{\textup{#3(#1)#4--#5(#2)#6}}%
{ and \textup{#3(#1)#4--#5(#2)#6}}{, \textup{#3(#1)#4--#5(#2)#6}}{, and \textup{#3(#1)#4--#5(#2)#6}}

\Crefformat{equation}{Equation~\textup{#2(#1)#3}}
\Crefrangeformat{equation}{Equations~\textup{#3(#1)#4--#5(#2)#6}}
\Crefmultiformat{equation}{Equations~\textup{#2(#1)#3}}{ and \textup{#2(#1)#3}}
{, \textup{#2(#1)#3}}{, and \textup{#2(#1)#3}}
\Crefrangemultiformat{equation}{Equations~\textup{#3(#1)#4--#5(#2)#6}}%
{ and \textup{#3(#1)#4--#5(#2)#6}}{, \textup{#3(#1)#4--#5(#2)#6}}{, and \textup{#3(#1)#4--#5(#2)#6}}

\crefname{section}{Section}{Sections}

\usepackage{tikz}
\tikzset{normalnode/.style={circle, draw, fill=black, inner sep=0, minimum width=1.5mm}}

\usetikzlibrary{fit}
\usetikzlibrary{shapes,arrows,shadows,positioning}
\usetikzlibrary{backgrounds}

\newcommand{\ie}{i.e.}

\newcommand{\whp}{w.h.p.\xspace}

\newcommand{\Nn}{\mathbb{N}}

\newcommand{\ber}[1]{\operatorname{Ber}\!\left(#1\right)}

\newcommand{\geo}[1]{\operatorname{Geo}\!\left(#1\right)}

\newcommand{\bin}[2]{\operatorname{Bin}\!\left( #1,#2 \right)}

\newtheorem{theorem}{Theorem}[section]

\newtheorem{proposition}[theorem]{Proposition}

\newtheorem{problem}[theorem]{Problem}
\newtheorem{lemma}[theorem]{Lemma}
\newtheorem{observation}[theorem]{Observation}

\theoremstyle{definition}

\theoremstyle{remark}

\newtheorem{claim}[theorem]{Claim}

\newenvironment{pocd}[1]{\begin{proof}[Proof of {C}laim #1.]}{\end{proof}}

\newcommand{\ind}{\sqsubseteq}

\renewcommand{\le}{\leqslant}

\renewcommand{\leq}{\leqslant}
\renewcommand{\geq}{\geqslant}

\renewcommand{\Pr}[1]{\mathbb{P}\left(#1\right)}
\newcommand{\Ex}[1]{\mathbb{E} \left(#1\right)}

\newcommand{\comment}[1]{}

\renewcommand{\epsilon}{\varepsilon}

\newcommand{\fun}{\mathtt{fun}}
\newcommand{\vfun}{\mathrm{fun}}
\newcommand{\rep}{\mathtt{rep}}
\newcommand{\dgn}{\mathtt{dgn}}

\title{Functionality of Random Graphs}

\date{}

\author{John Sylvester\thanks{Department of Computer Science, University of Liverpool, UK, \texttt{john.sylvester@liverpool.ac.uk}, \orcidlink{0000-0002-6543-2934}
}
	\and
	Viktor Zamaraev\thanks{Department of Computer Science, University of Liverpool, UK, \texttt{viktor.zamaraev@liverpool.ac.uk}, \orcidlink{0000-0001-5755-4141}
	}
	\and Maksim Zhukovskii\thanks{School of Computer Science, University of Sheffield, UK, \texttt{m.zhukovskii@sheffield.ac.uk}, \orcidlink{0000-0001-8763-9533}
	}}

\begin{document}

	\maketitle
	\begin{abstract}
		The functionality of a graph $G$ is the minimum number $k$ such that in every induced subgraph
		of $G$ there exists a vertex whose neighbourhood is uniquely determined by the neighborhoods of at most $k$ other vertices in the subgraph.
		The functionality parameter was introduced in the context of adjacency labeling schemes, and it generalises a number of classical and recent graph parameters including degeneracy, twin-width, and symmetric difference. We establish the functionality of a random graph $G(n,p)$ up to a constant factor for every value of~$p$. 
	\end{abstract}
	
	\section{Introduction}
	\label{sec:intro}
	
	A common approach in algorithmic and structural graph theory to make a problem tractable is by considering the problem on restricted classes of graphs, for example planar graphs or interval graphs.
	One way to obtain a restricted graph class is to fix a graph parameter and consider the graphs where the parameter is bounded from above by a fixed constant, for example graphs of bounded degree or bounded treewidth.
	
	Functionality is a graph parameter that was introduced in the context of adjacency labeling schemes \cite{ACLZ15, AAL21}.
	Speaking informally, the functionality of a graph $G$, denoted $\fun(G)$, is the minimum number $k$ such that every induced subgraph of $G$ contains a vertex whose neighbourhood is uniquely determined by the neighbourhoods of at most $k$ other vertices in the subgraph. This parameter generalises maximum degree, degeneracy, twin-width, and symmetric difference in the sense that if any of these parameters is bounded in a class of graphs, then functionality is also bounded in the class. Despite its generality, functionality is still a relatively restrictive graph parameter as any graph class of bounded functionality is at most factorial, i.e., contains $2^{O(n \log n)}$ graphs on $n$ vertices \cite{ACLZ15}, which is significantly smaller than the number $2^{\binom{n}{2}}$ of all $n$-vertex graphs. However, the structure of graph classes with bounded functionality is far from being understood.
	
	It is known that the functionality of any $n$-vertex graph is at most $O(\sqrt{n \log n})$ and  the incidence graphs of projective planes have functionality $\Omega(\sqrt{n})$ \cite{DFOPSA23}.
	A number of studies investigated functionality of graphs from some special classes. It is known that unit interval graphs have functionality at most $2$, line graphs have functionality at most $6$, permutation graphs have functionality at most $8$ \cite{AAL21}, and interval graphs have functionality at most $8$ \cite{DLMSZ23}. On the other hand there are examples of graph classes that are at most factorial, but have unbounded functionality, e.g., hypercube graphs \cite{AAL21} and box intersection graphs in $\mathbb{R}^d$ for any $d \geq 3$ \cite{DLMSZ23}  have unbounded functionality.
	It was shown in \cite{DFOPSA23} that the functionality of a random graph $G(n,p)$ is $\Omega(\log n)$ with high probability\footnote{With probability approaching 1 as $n\to\infty$.} (\whp) for any $p \in \left(\frac{3\log n}{n}, 1-\frac{3 \log n}{n}\right)$, however nothing better than the general upper bound of $O(\sqrt{n\log n})$ was known for the functionality of random graphs. 
	
	\subsection{Our Contribution \& Discussion}
	In this paper, we determine functionality of $G(n,p)$ up to a constant factor for the entire range of $p \in [0,1]$. To avoid complicating the statement of \cref{thm:main} we state this only for $p=\Omega(1/n)$, since for $p=o(1/n)$ \whp $G(n,p)$ is acyclic and thus its functionality is at most $1$ (in particular, since functionality does not exceed the degeneracy). We note also that since functionality is invariant under complementation, i.e.~$\fun(G) =\fun(\overline{G})$ for any graph $G$, it suffices to consider $p\in [0,1/2]$. 
	
	\begin{theorem}[Main]\label{thm:main}
		Let $G_n \sim G(n,p)$, where $\Omega(1/n)= p:=p(n)\leq 1/2$. Then for any $w := w(n)\geq 1$,  \whp
		\begin{equation*}
			\fun(G_n) =
			\begin{cases}
				\Theta\left( \frac{np \ln(ew)}{\ln n} \right) & \text{ if } p = \sqrt{\frac{\ln n}{n w}}; \\
				\Theta\left( \frac{\ln (ew)}{p} \right) & \text{ if } p = \sqrt{\frac{w\ln n}{n}}.
			\end{cases}
		\end{equation*}
	\end{theorem}
	
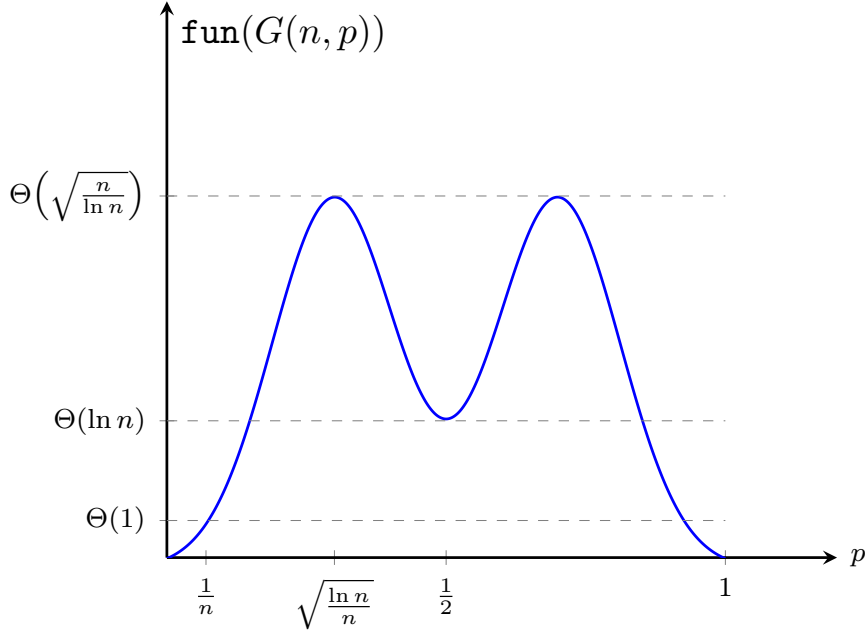
\begin{figure}[ht]
		\centering


	\scalebox{1.3}{
	\begin{tikzpicture}
		\begin{axis}[
			axis lines=middle,
			xlabel={$p$},
			ylabel={$\fun(G(n,p))$},
			xlabel style={font=\scriptsize},
			ymin=0, ymax=1.5,
			xmin=0, xmax=1.2,
			samples=200,
			domain=0:1,
			smooth,
			thick,
			xlabel style={at={(1,0)}, anchor=west}, 
			xtick={0.07, 0.3, 0.5,1}, 
			xticklabels={$\frac{1}{n}$,$\sqrt{\frac{\ln n}{n}}$,$\frac{1}{2}$,1}, 
			%
			ytick={0.1, 0.36887, 0.975}, 
			yticklabels={$\Theta(1)$,$\Theta(\ln n)$,$\Theta\Big(\sqrt{\frac{n}{\ln n}}\Big)$}, 
			tick label style={font=\scriptsize},
			]
			\addplot[
			blue,
			]
			{exp(-40*(x-0.3)^2) + exp(-40*(x-0.7)^2)-0.03};
			
			\pgfmathsetmacro{\fmin}{exp(-40*(0.5-0.3)^2) + exp(-40*(0.5-0.7)^2)-0.03}
			\pgfmathsetmacro{\fmax}{exp(0)}
			
			\addplot[
			dashed,
			gray,
			line width=0.3pt,
			domain=0:1
			]
			{0.36887};
			
			\addplot[
			dashed,
			gray,
			line width=0.3pt,
			domain=0:1
			]
			{0.975};
			
			\addplot[
			dashed,
			gray,
			line width=0.3pt,
			domain=0:1
			]
			{0.1}; 
			
		\end{axis}
	\end{tikzpicture}
}
		\caption{Qualitative behaviour of the functionality of $G(n,p)$}
		\label{fig:fun-plot}
	\end{figure}

	From \cref{thm:main}, one can see that around the point $p^* = \sqrt{\frac{\ln n}{n}}$ the functionality of $G(n,p)$ is maximised and becomes $\Theta(\sqrt{\frac{n}{\ln n}})$.
	In the interval $p \in [0, p^*]$ the functionality grows from $0$ to $\Theta(\sqrt{\frac{n}{\ln n}})$, and in the interval $p \in [p^*, 1/2]$ it decreases from $\Theta(\sqrt{\frac{n}{\ln n}})$ to $\Theta(\ln n)$. The behaviour on the interval $p \in [1/2, 1]$ mirrors that on the interval $p \in [0, 1/2]$ since, as mentioned above, $\fun(G) = \fun(\overline{G})$ holds for every graph $G$. \cref{fig:fun-plot} shows qualitative behaviour of the functionality of $G(n,p)$.
	
	As we have observed, the maximum functionality of $G(n,p)$ is achieved around $p^* = \sqrt{\frac{\ln n}{n}} = o(1)$  and not at $p = \Theta(1)$  where all other hereditary parameters mentioned above (namely, maximum degree, degeneracy, twin-width, symmetric difference) are maximised.
	In fact, the functionality of $G(n,p)$ drops  from $\Theta(\sqrt{\frac{n}{\ln n}})$ to $\Theta(\ln n)$ as $p$ increases from $p^*$ to $\Theta(1)$.
	This behaviour might at first seem mysterious, but it should hopefully be clear from the proofs of our bounds, outlined in \cref{sec:proofstrat}. In very rough terms, what is going is that for $p$ below $p^*$, if one wishes to encode the neighbourhood $N(v)$ of a vertex $v$ in a random graph using adjacencies with some set of vertices $S\subseteq V$, then one cannot do significantly better than just taking $S=\{v\}\cup N(v)$. However, for $p$ above $p^*$, one can find a set $S$ with $|S|\ll |N(v)|$ where for each $x\in N(v)$ and $y\notin N(v)$, the adjacencies of $x$ and $y$ with the set $S$ are different. 
	
	 Second, random graphs $G(n,p)$ almost achieve maximum functionality among $n$-vertex graphs, up to a log-factor. Indeed, the maximum functionality of a random graph is  $\Theta(\sqrt{\frac{n}{\ln n}})$ \whp, and the best known general upper bound on functionality \cite{DFOPSA23} is $O(\sqrt{n\ln n})$. It is interesting to note that incidence graphs of projective planes have functionality $\Theta(\sqrt{n})$ \cite{DFOPSA23}, and thus, in contrast with many other width measures such as treewidth, rank-width \cite{LeeLO12}, and twin-width \cite{ACHDO24} the functionality is not maximised (up to a constant) by some choice of parameters for the random graph. 
	
 	We discuss how this result is proved in the next subsection. However, here we state the following result on dominating sets in random bipartite graphs as we believe it may be of independent interest. We used a stronger version of this theorem (\cref{lem:small_t_dom_exist}) in our proofs.

	\begin{restatable}{theorem}{simpleDomLemma}\label{lem:small_t_dom_simple}
		Let $\alpha, \beta>0$ be fixed constants,  $n$ be sufficiently large, and $f(n)$ be any function satisfying $f(n)\rightarrow \infty$ as $n\rightarrow \infty$. Let  $\ln^2 n\leq a\leq b\leq n$ and $p\in(0,1]$ be functions of $n$ which satisfy $ap \leq   \alpha\cdot \ln n $ and $bp >n^{\beta}$.  Then, for $G\sim G(A,B,p)$, where $|A|=a $ and  $|B|=b$, we have	
		\[	
			\Pr{ \exists \text{ a dominating set $D\subset B$ of }A\text{ with }\left. |D|\leq \frac{f(n)}{p} \, \right| \, \forall a\in A,\; N(a)\neq \emptyset } \geq  1- \exp\left[- b^{1-o(1)}\right]. 
		\] 
	\end{restatable}
	
We note that if one is willing to accept a less strong probability guarantee then one can exchange $f(n)$ for a suitably large constant $C>0$ to obtain a simpler bound from \cref{lem:small_t_dom_exist}, such as
\[
	\exists C \text{ s.t. } \Pr{\exists \text{ a dominating set $D\subset B$ of }A\text{ with }\left. |D|\leq \frac{C}{p} \, \right| \, \forall a\in A,\; N(a)\neq \emptyset } 
	\geq  1- \exp\left[- b^{0.9}\right]. 
\]

We further observe that \cref{lem:small_t_dom_simple} is optimal in the following sense. First, \whp a minimum dominating set of $A$ in $B$ has size $\Omega(1/p)$ whenever $ap=\Theta(\ln n)$ (see \cref{prop:dom-lower-bound}). Second, the probability that there exists a dominating set of size $O(1/p)$ is at most $1-\exp\left[-O(b)\right]$ as soon as $ap=\Omega(1)$ (see \cref{prop:dom-upper-bound}).
 
We are not aware of any such result for dominating sets in random bipartite graphs, but for related results on dominating sets in random graphs see \cite{NikoletseasS93,bohman2024,GlebovLS15,Weber81,WielandG01}.

	\subsection{Proof Strategies}\label{sec:proofstrat}
	
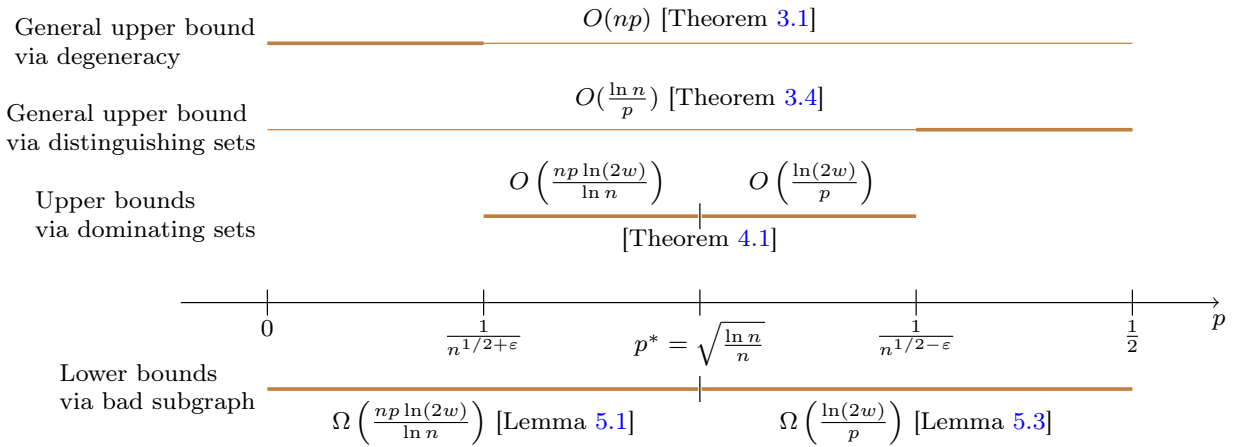
\begin{figure}[ht]
		\centering

	\scalebox{1.15}{
	\begin{tikzpicture}
		
		\draw[thin,brown] (0,4) -- (10,4);
		\draw[very thick,brown] (0,4) -- (2.5,4);
		\node[left, font=\scriptsize, align=left] at (0,4) {General upper bound \\ via degeneracy}; 
		\node[above, font=\scriptsize, align=left] at (5,4) {$O(np)$ [\cref{th:general-via-dgn}]}; 
		
		\draw[thin,brown] (0,3) -- (10,3);
		\draw[very thick,brown] (7.5,3) -- (10,3);
		\node[left, font=\scriptsize, align=left] at (0,3) {General upper bound \\ via distinguishing sets}; 
		\node[above, font=\scriptsize, align=left] at (5,3) {$O(\frac{\ln n}{p})$ [\cref{th:general-via-distinguishing-set}]}; 
		
		\draw[very thick,brown] (2.5,2) -- (4.98,2);
		\draw[very thick,brown] (5.02,2) -- (7.5,2);
		\draw[thin] (5,1.85) -- (5,2.15); 
		\node[left, font=\scriptsize, align=left] at (0,2) {Upper bounds \\ via dominating sets}; 
		
		\draw[thin,->] (-1,1) -- (11,1);
		\draw[thin] (0,0.85) -- (0,1.15) node[pos=0.25, below, font=\scriptsize] {0}; 
		\draw[thin] (2.5,0.85) -- (2.5,1.15) node[pos=0.25, below, font=\scriptsize] {$\frac{1}{n^{1/2+\epsilon}}$}; 
		\draw[thin] (5,0.85) -- (5,1.15) node[pos=0.25, below, font=\scriptsize] {$p^* = \sqrt{\frac{\ln n}{n}}$}; 
		\draw[thin] (7.5,0.85) -- (7.5,1.15) node[pos=0.25, below, font=\scriptsize] {$\frac{1}{n^{1/2-\epsilon}}$}; 
		\draw[thin] (10,0.85) -- (10,1.15) node[pos=0.25, below, font=\scriptsize] {$\frac{1}{2}$}; 
		\node[below, font=\scriptsize] at (11,1) {$p$}; 
		\node[above, font=\scriptsize] at (3.7,2) {$O\left( \frac{np \ln(2w)}{\ln n}\right)$}; 
		\node[below, font=\scriptsize, align=left] at (5,2) {[\cref{th:fun_upper2}]};
		\node[above, font=\scriptsize] at (6.3,2) {$O\left( \frac{\ln (2w)}{p} \right)$}; 
		
		\node[left, font=\scriptsize, align=left] at (0,0) {Lower bounds \\ via bad subgraph}; 
		\draw[thin] (5,-0.15) -- (5,0.15); 
		\draw[very thick,brown] (0,0) -- (4.98,0);
		\draw[very thick,brown] (5.02,0) -- (10,0);
		\node[below, font=\scriptsize] at (2.5,0) {$\Omega\left( \frac{np \ln(2w)}{\ln n}\right)$ [\cref{lem:lower-p-less-p-star}]}; 
		\node[below, font=\scriptsize] at (7.5,0) {$\Omega\left( \frac{\ln (2w)}{p} \right)$ [\cref{lem:lower-p-more-p-star}]}; 
		
	\end{tikzpicture}
}
		\caption{Upper and lower bounds on the functionality of $G(n,p)$. The \textcolor{brown}{brown} segments show the intervals on which the corresponding bounds hold. The \textcolor{brown}{\textbf{thick brown}} segments show the intervals where the bounds are order optimal.}
		\label{fig:bounds}
	\end{figure}

	Obtaining order optimal bounds on functionality requires different arguments for different regimes of $p$.
	The summary of the bounds is depicted in \cref{fig:bounds}.
	We now proceed with an outline of the proof strategies for different bounds.
	We consider $p$ in the interval $[0,1/2]$ as the situation in the interval $[1/2,1]$ is symmetric.
	
	\paragraph{General upper bounds (small and large values of $p$).} 
	We first obtain two general upper bounds that hold for the entire interval $[0, 1/2]$. The first bound is order optimal on the interval $[0, n^{-1/2-\epsilon}]$, and the second bound is order optimal on $[n^{-1/2+\epsilon}, 1/2]$, for any fixed $\varepsilon>0$.
	
	The first bound is based on the observation that any vertex is a function of its neighbourhood, and therefore $\fun(G)$ does not exceed the degeneracy of $G$, which is the maximum of the minimum degree of a subgraph of $G$. Thus, 
	we obtain a bound of $O(np)$ on the functionality of $G(n,p)$ from the same bound on its degeneracy (\cref{lem:dgn}).
	
	The second bound is based on the observation (\cref{obs:fun2}) that if $S \subseteq V(G)$ is
	a \emph{distinguishing set} in $G$, i.e., a set such that all vertices in $V(G) \setminus S$ have pairwise different neighbourhoods in $S$, then any vertex $y \in V(G) \setminus S$ is a function of $S$. 
	More specifically, in \cref{th:general-via-distinguishing-set}, we show that in $G_n \sim G(n,p)$ \whp~for every subset $S' \subseteq V(G_n)$ of size $\lfloor \frac{3\ln n}{p} \rfloor$, the set $R \subseteq V(G_n) \setminus S'$ of vertices with a non-unique neighbourhood in $S'$ has at most $\frac{7\ln n}{p}$ elements.
	This is enough to conclude (see \cref{lem:repitionbound}) that 
	in any induced subgraph $H$ of $G_n$ there exists a set $S$
	of at most $\frac{10\ln n}{p}$ vertices such that any vertex in $V(H) \setminus S$ is a function of $S$ and thus \whp~$\fun(G_n) \leq \frac{10\ln n}{p}$.

		\paragraph{Upper bounds via dominating sets ($p$ around $p^*$).} 
			In the middle interval $[n^{-1/2-\epsilon}, n^{-1/2+\epsilon}]$, the above two general approaches (degeneracy and the maximum size of a minimum distinguishing set over all induced subgraphs) do not give order optimal bounds on functionality, and we need to use a more refined argument to find in every induced subgraph a vertex that is a function of some set of vertices of the appropriate size, for the sake of this discussion let us denote this size $d:=d(p)$.

To bound the functionality of a vertex $y$ one can find a set $S$ such that any pair of vertices that have different adjacencies with $y$ also have different neighbourhoods in $S$ (i.e., the two vertices are distinguished by $S$). We use this idea to obtain the upper bound on the functionality for $p$ around $p^*$. Namely, we build such a set $S$ of order $d$ by uniting several dominating sets of the neighbourhood of a vertex $y$. We show that such a set exists with large enough probability that it is also present in all induced subgraphs. The main steps are:

\begin{enumerate}[(i)] 
	\item\label{step1} We prove a general  bound on  dominating sets in  random bipartite graphs (\cref{lem:small_t_dom_exist}).   
	\item\label{step2} We show that given at least five distinct sets of order $d$ which dominate the neighbourhood of a vertex $y$, \whp~there are only order $d$ many non-neighbours of $y$ that cannot be distinguished from a neighbour of $y$ by the union of these sets  (\cref{lem:induced-domination}). Thus, if we take these five dominating sets together with the vertices which are not distinguished by them, then we have a set $S$ of order $d$ such that $y$ is a function $S$. 
	\item \label{step3} We then apply the general result for domination in random bipartite graphs from \eqref{step1} to find \whp~a dominating set of order $d$ of the neighbourhood of a given vertex $y$ (\cref{cl:from_AB_reduction}). 
	\item\label{step4} We then use this to show that there exists a vertex whose neighbourhood has at least five distinct dominating sets of order $d$  (\cref{lem:clm1}) with even higher probability. 
	We achieve this higher probability of the event, over what would be guaranteed by \eqref{step3} for a fixed vertex, by using the fact that most neighbourhoods in our random graph do not overlap much.
	We do this probability boosting to overcome a large union bound in \eqref{step5}. 
	\item\label{step5} We then wrap up the proof in \cref{sc:TH_upper_inter_proof}, dealing with `large' and `small' induced subgraphs separately. 
	Since the probability guarantee by \eqref{step4} is so strong, we can take a union bound over all `large' induced subgraphs to show that \whp~each of them contains a vertex $y$ whose neighbourhood is dominated by at least five distinct sets of order $d$, which, by \eqref{step2}, can be converted into a witness of $y$ having functionality of order $d$. For `small' induced subgraphs we use an elementary bound on the degeneracy to bound the functionality. 
\end{enumerate}

	\paragraph{Lower bounds.}
	We obtain our lower bounds by showing the existence \whp~of a subgraph $H$ in $G_n \sim G(n,p)$ such that no vertex $y$ in $H$ is a function of a ``small'' set of vertices.
	To do so, we show that \whp~there exists $H$ such that for every ``small'' set $S \subseteq V(H)$ and every vertex $y \in V(H)$ outside $S$ there exists a neighbour $z_1$ and a non-neighbour $z_2$ of $y$ in $H$ that have no neighbours in $S$. This means that $S$ cannot distinguish $z_1$ and $z_2$ and therefore cannot be used to represent different adjacencies between $y$ and $z_1$, and $y$ and $z_2$, i.e., $y$ cannot be a function of~$S$.
	This strategy requires slightly different implementations in the intervals $[0, p^*]$ (\cref{sec:lowersmallp}) and
	$[p^*,1]$ (\cref{sec:lowerlargep}). Aside from technical differences arising from differing rates of decay for $p$, the main distinction is that for small $p$ the graph $H$ is a proper subgraph of $G_n$, whereas for larger $p$ we can take $H=G_n$.

	\subsection{Related Graph Parameters}
	
	As was already mentioned, functionality generalises a number of graph parameters, in the sense that functionality is upper bounded by a function of these parameters, and, thus, if a parameter is bounded in a class of graphs, then functionality is bounded too.
	Among these parameters are maximum degree, degeneracy, symmetric difference, and twin-width.
	All of these parameters are hereditary, i.e., their value on induced subgraphs cannot exceed the value on the graph.
	The maximum degree (and even the degree distribution) of a random graph $G(n,p)$ is well-understood in the entire range $p \in [0,1]$ (see e.g.~\cite[Chapter 3]{RandomGraphs}).
	The behaviour of the recently introduced parameter twin-width of random graphs has been studied in \cite{AHKO22,ACHDO24,HNST23}. 
	
	Two more parameters, similar in nature to functionality, that were investigated for random graphs are metric dimension \cite{BMP13} and the minimum size of identifying sets \cite{FMMRS07}. 
	The metric dimension of a graph $G$ is the minimum size of a vertex set $S$ of $G$ such that the vertices outside $S$ can be distinguished by their distances to the vertices in $S$. Similarly, the minimum size of an identifying set is defined as the size  of a smallest vertex set $S$ such that the vertices outside $S$ are distinguished by their adjacencies to the vertices in $S$.
	The main conceptual difference between these notions and functionality is that in the definition of the latter a ``distinguishing'' set must exist in every induced subgraph, rather than only in the graph itself.
	This hereditary property of functionality poses new challenges, both quantitative and qualitative, that do not arise in the analysis of metric dimension and identifying sets.
	In particular, achieving a tight upper bound necessitates overcoming the union bound over $2^n$ induced subgraphs, which demands at least exponentially small probability bounds for rare events. 
	This makes the proofs of the upper bounds especially intricate and requires establishing new tight domination results in bipartite random graphs (\cref{lem:small_t_dom_exist}).

	\section{Preliminaries}
	\label{sec:pre}
	
	For a positive integer $k$, we denote by $[k]$ the set $\{1,\ldots, k\}$.
	
	\paragraph{Graphs.}
	We consider only finite, undirected graphs, without loops and multiple edges.
	The vertex set and the edge set of a graph $G$ are denoted $V(G)$ and $E(G)$, respectively.
	The {\it neighbourhood} of a vertex $x\in V(G)$, denoted $N(x)$, is the set of vertices of $G$ adjacent to $x$, and the {\it degree} of $x$, denoted  $\deg(x)$, is the size of its neighbourhood.
	The minimum vertex degree and the maximum vertex degree in $G$ are denoted by $\delta(G)$ and $\Delta(G)$, respectively.
	We denote by $N[x]$ the {\it closed neighbourhood} of $x$, that is, the set $N(x)\cup \{x\}$. For a set $S \subseteq V(G)$ we write $N_S(x)$ to denote the neighbourhood of $x$ in $S$, \ie, $N_S(x) = N(x) \cap S$. 
	For a set $S \subseteq V(G)$ we denote by $G[S]$ the subgraph of $G$ induced by $S$. We write $H \ind G$ to denote the fact that $H$ is an induced subgraph $G$.
	
	Given a set $S \subseteq V(G)$ and a number $t \in \mathbb{N}$, a set $D \subseteq V(G) \setminus S$ is said to \emph{$t$-dominate} $S$ if every vertex in $S$ has at least $t$ neighbours in $D$; in which case we also say that $D$ is a \emph{$t$-dominating set} for $S$. A set $D \subseteq V$ is a \emph{$t$-dominating set} in $G$ if it $t$-dominates $V \setminus D$. If $t=1$, then we simply say ``dominates'' and ``dominating set'' instead of ``$t$-dominates'' and ``$t$-dominating set'', respectively.
	
 Let $G(n,p)$ denote the distribution on all graphs on $[n]$ where each edge is included independently with probability $p:=p(n)\in[0,1]$, i.e., for every graph $G$ on $[n]$ and $G_n\sim G(n,p)$, 
 $$
 \mathbb{P}(G_n=G)=p^{|E(G)|}(1-p)^{{n\choose 2}-|E(G)|}.
 $$
Similarly for two disjoint sets $A,B$ and $p:=p(n)\in[0,1]$, we let $G(A,B,p)$ denote the random bipartite graph where each edge with one end point in $A$ and the other in $B$ is included independently with probability $p$. 
	
	\paragraph{Functionality.}
	Let $G$ be a graph and let $A=A_G$ be the adjacency matrix of $G$.
	We say that a vertex $y\in V(G)$ is a {\it function of (a set of) vertices} $S=\{x_1,x_2,\ldots,x_k\} \subseteq V(G) \setminus\{y\}$ if there is a Boolean function $f$ of $k$ variables such that for every vertex $z\in V(G)$ different from $y,x_1,x_2,\ldots,x_k$ we have $A(y,z)=f(A(z,x_1),\ldots,A(z,x_k))$.
	Alternatively, $y$ is a function of $S$ if for every $z \in V(G) \setminus (S \cup \{y\})$ the adjacency between $y$ and $z$ in $G$ depends only on $N(z) \cap S$. 
	
	The following observations are simple implications of definitions.
	\begin{observation}\label{obs:fun1}
		Let $G=(V,E)$ be a graph, $y \in V$, and $S \subseteq V \setminus \{y\}$.
		Then $y$ is a function of $S$ if and only if there does not exist a pair of vertices $z_1,z_2 \notin S \cup \{y\}$ such that $N_S(z_1) = N_S(z_2)$, and $z_1$ and $z_2$ have different adjacencies with $y$.
	\end{observation} 
	
	For a graph $G=(V,E)$, a set of vertices $S \subseteq V$ is called \emph{distinguishing set} in $G$ if no two vertices outside $S$ have the same neighbourhood in $S$, i.e., $N(v) \cap S \neq N(w) \cap S$ for any two distinct vertices $v,w \in V \setminus S$. 
	
	\begin{observation}\label{obs:fun2}
		Let $G$ be a graph and $S \subseteq V(G)$ be a distinguishing set in $G$.
		Then any vertex $y \in V(G) \setminus S$ is a function of $S$.
	\end{observation}
	
	We observe that every vertex is a function of some other vertices.
	In particular, every vertex is a function of its neighbours, for the function $f\equiv 0$, and a function of its non-neighbours,  for the function $f\equiv 1$.
	The minimum $k$ such that $y$ is a function of $k$ other vertices is the {\it functionality of $y$} and is denoted $\vfun_G(y)$, or simply $\vfun(y)$ if the graph is clear from the context. The \emph{functionality of $G$} is denoted and defined as follows:
	\begin{equation}
		\fun(G)=\max_{H \ind G} \min_{y\in V(H)} \vfun(y),
		\label{eq:fun-definition}
	\end{equation}
	If in this definition we replace $\vfun(y)$ with $\deg(y)$, then we obtain the definition of \emph{degeneracy}:
	
	$$
		\dgn(G)=\max_{H \ind G} \min_{y\in V(H)} \deg(y).
	$$
	
	From the two definitions and the observation that $\vfun(y) \le \deg(y)$, we immediately obtain
	\begin{observation}\label{obs:fun-dgn}
		For any graph $G$ it holds $\fun(G) \leq \dgn(G)$.
	\end{observation}
	
	\paragraph{Probability.} 
	
	For random variables $Y,Z$ we say that $Y$ is stochastically dominated by $Z$, denoted by $Y \preceq Z$, if $ \mathbb{P}( Y \geq x) \leq \mathbb{P}(Z \geq x )$ for all real $x$. We use $\mathrm{Bin}(n,p)$ for the binomial distribution with $n\geq 1$ trials, each with success probability $p\in [0,1]$. We now state some standard Chernoff bounds for the binomial distribution.

	\begin{lemma}[{\cite[Theorem 2.1]{JansonBook}}]	\label{lm:Chernoff_JLR}
		Let $n\in\mathbb{N}$, $p\in[0,1]$, $X\sim\mathrm{Bin}(n,p)$, $t>0$. Then \[\Pr{X<np -t}\leq \exp\left[-\frac{t^2}{2np}\right] \qquad \text{and}\qquad \Pr{X>np +t}\leq \exp\left[-\frac{t^2}{2(np +t/3)}\right].\] 
		
	\end{lemma}

	\begin{lemma}[{\cite[Theorem A.1.12]{AlonSpencer}}]\label{lm:Chernoff_AS}
		Let $n\in\mathbb{N}$, $p\in[0,1]$, $X\sim\mathrm{Bin}(n,p)$, $\lambda>0$. Then \[\Pr{X>\lambda np}\leq (e/\lambda)^{\lambda np}.\]
	\end{lemma}
	Let $\geo{p}$ denote the \textit{geometric distribution} with success probability $p\in (0,1)$. That is, if a random variable $X\sim \geo{p}$ then $\Pr{X=k} = (1-p)^{k-1}p  $ for any integer $k\geq 1$. We will also need the following Chernoff bound for sums of geometric random variables. 
	\begin{lemma}[{\cite[Theorem 2.1]{JansonTail}}] \label{lem:jansontail}
		For any $n\geq 1$ and $p_1,\dots  , p_n \in (0,1)$, let $X=\sum_{i=1}^nX_i$,
		where $X_i \sim \geo{p_i}$, $i \in [n]$, are  independent geometric random variables. Let $p_*=\min_{i\in [n]}p_i$ and $\mu=\Ex{X}=\sum_{i=1}^n\frac{1}{p_i}$. 
		Then for any $\lambda\geq 1$,  
		\[
			\Pr{X \geq \lambda \mu  } \leq \exp\left[ -p_*\cdot \mu\cdot  \left(\lambda -1- \ln \lambda\right)\right].
		\]
	\end{lemma}

	\paragraph{Properties of random graphs.}	
	We require the following two standard technical results, which we prove here for completeness. 

	\begin{lemma}\label{lem:dgn}
		Let $p:=p(n) \in [0,1]$ and $G_n\sim G(n,p)$. Then \whp~$\dgn(G_n)\leq \max\{10np, 2\}$.
	\end{lemma}
	
	\begin{proof}
		If $p < 0.5/n$, then \whp~every connected component of $G_n$ contains at most one cycle \cite[Theorem 5.5]{JansonBook}, which implies $\dgn(G_n) \leq 2$. Thus, we now assume that $p \geq 0.5/n$ and show that $\dgn(G_n) \leq 10np$.
		
		It is sufficient to prove that \whp, for every $H\ind G_n$, $\frac{|E(H)|}{|V(H)|}\leq 5np$. Let us fix $k\in[n]$ and consider any $k$-element subset $V\subset [n]$. Then the number of edges induced by $V$ is $\xi\sim\mathrm{Bin}({k\choose 2},p)$. By Lemma~\ref{lm:Chernoff_AS},
		$$
		\Pr{\xi>5npk}\leq\left(\frac{e{k\choose 2}p}{5npk}\right)^{5npk}\leq \exp\left[-5npk\ln\frac{9n}{ke}\right].
		$$
		By the union bound, probability that there exists a subset $V\subset [n]$ of size $k$ that induces more than $5npk$ edges is at most
		$$
		{n\choose k}\Pr{\xi>5npk}\leq \exp\left[k\left(\ln\frac{ne}{k}-5np\ln\frac{9n}{ke}\right)\right]
		\leq \exp\left[k\left(\ln\frac{ne}{k}-2\ln\frac{9n}{ke}\right)\right]
		\leq \exp\left[-k\ln\frac{9n}{ke}\right].
		$$
		Finally, by the union bound,
		$$
		\Pr{\exists V\subset[n]\,:\,\frac{|E(H)|}{|V(H)|}>5np}\leq\sum_{k=1}^n \exp\left[-k\ln\frac{9n}{ke}\right]
		\leq\sum_{k=1}^{\lfloor\sqrt{n}\rfloor} \exp\left[-\frac{1}{2}k\ln n \right]+\sum_{k=\lceil\sqrt{n}\rceil}^{n} \exp\left[-k\right],
		$$which is $o(1)$, as claimed. 
	\end{proof}

	\begin{lemma}
		Let $p\gg\frac{\ln n}{n}$ and $G_n\sim G(n,p)$. Then \whp~every vertex has degree $np(1+o(1))$ and every two vertices have at most $\max\{2\ln n,10np^2\}$ common neighbours in $G_n$.
		\label{cl:degrees_common_degrees}
	\end{lemma}
	
	\begin{proof}
		Let $w=w(n)\to\infty$ as $n\to\infty$ be such that $p/w\gg\frac{\ln n}{n}$. By the Chernoff bound (\cref{lm:Chernoff_JLR}), a fixed vertex has degree which is at least $\frac{np}{w}$-far from $np$ with probability at most $\exp[-\Omega(np/w)]$. Thus, by the union bound, \whp~every vertex has degree in the interval $[np-np/w,np+np/w]$. Observe that, for any fixed pair of vertices, the number of vertices in their common neighbourhood is stochastically dominated by a random variable with distribution $\bin{n}{p^2}$. Thus, by the Chernoff bound, if $np^2\geq \frac{1}{5}\ln n$, then a fixed pair of vertices have at least $10np^2$ common neighbours with probability at most $\exp\big[-\frac{(9np^2)^2}{2(np^2 + 3np^2)}\big]=o(n^{-2})$. On the other hand, if   $np^2<\frac{1}{5} \ln n$, then, again by the Chernoff bound, a fixed pair of vertices has at least $2\ln n$ common neighbours with probability at most $\exp\big[-\frac{(2\ln n)^2}{2(\frac{1}{5} \ln n + \frac{2}{3} \ln n)}\big]=o(n^{-2})$. The result thus follows from the union bound.
	\end{proof}

	\section{General Upper Bounds}\label{sec:general-bounds}
	In this section we establish two general upper bounds. The first one is based on degeneracy, and the second one is based on distinguishing sets.
	
	\paragraph{Upper bound via degeneracy.}
	The bound in \cref{th:general-via-dgn} below follows immediately from the fact that degeneracy upper bounds functionality  (\cref{obs:fun-dgn}) and the upper bound on degeneracy of $G(n,p)$ from \cref{lem:dgn}.
	
	\begin{theorem}\label{th:general-via-dgn}
		Let $p:=p(n) \in [0,1]$ and $G_n \sim G(n,p)$. Then \whp~$\fun(G_n) \leq \max\{10np, 2\}$.
	\end{theorem}

	\paragraph{Upper bound via distinguishing set.}
	For this bound we use the notion of a distinguishing set and the fact that any vertex outside a distinguishing set is a function of this set (\cref{obs:fun2}). In order to use this fact to upper bound the functionality, we show that any sufficiently large induced subgraph of $G(n,p)$ has a distinguishing set of the desired size.
	
	We proceed with some notation and auxiliary statements. For a graph $G=(V,E)$ and a set $U\subseteq V$ we define the \emph{repetition number of $U$ in $G$}, denoted $\rep_G(U)$, as the number of vertices outside $U$ that have the same neighbourhood in $U$ as some other vertex outside $U$, i.e.,
	\[
		\rep_G(U) := |\{ v \in V \setminus U :  \exists w \in V \setminus U \text{ s.t. } N_U(v) = N_U(w) \}|.
	\]
	
	The repetition number of the set $U$ shows how far $U$ is from being a distinguishing set in the sense that $U$ can be extended to a distinguishing set in $G$ by adding at most $\rep_G(U)$ vertices to it. 
	A simple, but important observation is that the repetition number of a set cannot increase in induced subgraphs.
	
	\begin{observation}\label{obs:repmnono} 
		For any graph $G$, $H \ind G$, and $U\subseteq V(H)$, we have $\rep_H(U) \leq \rep_G(U) $.   
	\end{observation}
	
	We use this observation to show that if every small set of vertices in a graph
	has small repetition number, then the graph has small functionality.
	
	\begin{lemma}\label{lem:repitionbound} 
		For any graph $G$ and $k,r \in \Nn$, if $\rep_G(U)\leq r$ for all $U\subseteq V(G)$ with $|U|= k$, then $\fun(G) \leq k+r$. 
	\end{lemma}
	\begin{proof}
		To prove the lemma, we will show that any induced subgraph of $G$ 
		has a vertex that is a function of at most $k+r$ other vertices.
		
		Consider any induced subgraph $H \ind G$ and note that we can assume that $V(H)\geq k+r+1$ or else any vertex of this subgraph is a function of at most $k+r$ other vertices. 
		Now, take any set $U\subseteq V(H)$ with $|U|=k$. By hypothesis we have $\rep_G(U)\leq r $, and thus $\rep_H(U)\leq r $ by \cref{obs:repmnono}. It follows that there is a set $R$ of at most $r$ vertices, such that all vertices in $V(H)\setminus (U \cup R)$ have distinct neighbourhoods in $U$, and therefore $U \cup R$ is a distinguishing set in $H$. 
		Hence, by \cref{obs:fun2}, any vertex in $V(H) \setminus (U \cup R)$ is a function of at most $|U| + |R| \leq k+ r$ vertices. 
	\end{proof}
	
	Using this result we establish our second general upper bound.  
	
	\begin{theorem}\label{th:general-via-distinguishing-set} 
		Let $p \in (0,1/2]$ and $G_n \sim G(n,p)$. Then \whp~$\fun\left(G_n\right) \leq \frac{10\log n}{p}$. 
	\end{theorem}
	\begin{proof} 
		Let $G_n=(V,E) \sim G(n,p)$.
		For $r,k\geq 0$, let the event $\mathcal{E}(k,r)$ be given by 
		\[
		\mathcal{E}(k,r) = \{\text{for all } U \subseteq V \text{ satisfying } |U|=k \text{ we have } \rep_{G_n}(U) \leq r  \}.
		\] 
		We prove the theorem by showing that for 
		$$
		k=\left\lfloor\frac{3\ln n}{p}\right\rfloor\quad\text{  and }\quad r=\frac{7\ln n}{p},
		$$
		the event $\mathcal{E}(k,r)$ holds \whp, which together with \cref{lem:repitionbound} implies  the result.
		
		Let $U \subseteq V$ be a set of vertices of size $k$.
		Fix an arbitrary ordering of vertices in $V \setminus U$, say $v_1,v_2, \ldots, v_{n-k}$, and for every $v_i \in V \setminus U$, define 
		\[
		R(G_n,U,v_i) := 
		\begin{cases}
			1 &\text{ if there exists some } j< i \text{ satisfying }  N_U(v_i) = N_U(v_j);\\
			0 &\text{ otherwise.} 
		\end{cases}
		\]
		Let $R(G_n,U) := \sum_{v_i\in V \setminus U} R(G_n,U,v_i)$ and observe that this quantity is independent of the ordering of $V\setminus U$.
		We claim that $\rep_{G_n}(U) \leq 2R(G_n,U)$. Indeed, let $C_1, C_2, \ldots, C_t \subseteq V \setminus U$ be the non-singleton equivalence classes of the vertices in $V \setminus U$, where two vertices are equivalent if and only if they have the same neighbourhood in $U$. Then, every such equivalence class $C_i$ contributes $|C_i|-1$ to $R(G_n,U)$ and $|C_i|$ to $\rep_{G_n}(U)$. Thus, we have
		\[
		R(G_n,U) = \sum_{i=1}^t (|C_i|-1)  \leq \rep_{G_n}(U) = \sum_{i=1}^t |C_i| \leq \sum_{i=1}^t 2(|C_i|-1) = 2R(G_n,U),
		\]
		where the second inequality holds, because $|C_i| \geq 2$ for every $i \in [t]$.
		
		Now, the probability that $v_i \in V\setminus U$ has a specific set of $j$ neighbours in $U$ is equal to $p^j(1-p)^{k-j}\leq (1-p)^k$, where the inequality holds because $p \leq 1/2$. As the worst case is when all previous vertices $v_1,\dots, v_{i-1}$ have had distinct neighbourhoods, for each 
		$i\leq k$,
		\[
		\Pr{R(G_n,U,v_i)=1}\leq i\cdot (1-p)^k\leq n \cdot e^{-pk} < \frac{2}{n^2}.
		\]
		Observe also that $R(G_n,U,v_i) \preceq \gamma_i$ where 
		$\gamma_i\sim \ber{  2/n^2}$ are i.i.d.\ Bernoulli random variables, thus $R(G,U)\preceq \bin{n}{2/n^2}$. 
		Using \cref{lm:Chernoff_AS} with $t= rn/2>0$, for any fixed $U$ we have 
		\begin{align*}
			\Pr{R(G_n,U)\geq r/2}
			&\leq \Pr{\bin{n}{2/n^2} \geq r/2} 
			\leq \left(\frac{4e}{rn}\right)^{r/2} 
			\stackrel{p\leq 1/2}\leq \left(\frac{e}{(3.5\ln n)n}\right)^{\frac{3.5\ln n}{p}} \leq   e^{- \frac{3.5(\ln n)^2}{p} }.
		\end{align*}

		We now observe that by the union bound 
		\begin{align*}
			\Pr{ \neg \mathcal{E}(k,r) }
			&\leq \binom{n}{k}\cdot \Pr{\rep_{G_n}(U)\geq r}  \leq
			\binom{n}{k}\cdot \Pr{R(G_n,U)\geq r/2} 
			\leq  n^{\frac{3\ln n}{p}}\cdot    e^{- \frac{3.5(\ln n)^2}{p} }  = o(1), 
		\end{align*} 
		which concludes the proof.
	\end{proof}

	\section{Tight Upper Bound for \texorpdfstring{$p$ around $p^*$}{p around p*}}
	
	Our aim in this section is to prove the following bound on functionality, which gives the right behaviour for $p$ close to the value $p^*$, i.e.\ near where the functionality of $G(n,p)$ is maximised. 
	\begin{theorem}\label{th:fun_upper2}
		Let $G_n \sim G(n,p)$, where $p\in [n^{-1/2-\varepsilon},n^{-1/2+\varepsilon}]$ and $\varepsilon >0$ is a sufficiently small constant. Then there exists a constant $c$ such that \whp
		\[
		\fun\left(G_n\right) \leq \left\{\begin{array}{cc} c\cdot \frac{np\ln(ew)}{\ln n}, &\text{ for } p=\sqrt{\frac{\ln n}{wn}},\, w=w(n)\geq 1; \\ c\cdot\frac{\ln(ew)}{p}, &\text{ for } p=\sqrt{\frac{w\ln n}{n}},\, w=w(n)\geq 1. \end{array}\right. 
		\] 
	\end{theorem}
	
	We prove Theorem~\ref{th:fun_upper2} by showing that, in all large enough induced subgraphs of $G(n,p)$ there exist a vertex $y$ and a set $D$ of an appropriate size that $t$-dominates most of the neighbourhood of $y$ such that there are only few non-neighbours of $y$ that have same neighbourhoods in $D$ as some neighbours of $y$. This proof heavily relies on an auxiliary result about dominating sets in random bipartite graphs with parts of different sizes stated below.
	
Recall that for two independent sets $A,B$ and $p\in[0,1]$, the random bipartite graph $G(A,B,p)$ denotes the distribution over bipartite graphs on $A,B$ where each edge $ab$ with $a\in A$ and $b\in B$, is present independently with probability $p$.

	\begin{restatable}{theorem}{DomLemma}\label{lem:small_t_dom_exist}
		Let $\alpha, \beta>0$, and $0<\delta <\beta$ be fixed constants, and $n$ be sufficiently large. Let  $a,b \in \mathbb{N}$ where $\ln^2 n\leq a\leq b\leq n$, and $h,p\in \mathbb{R}$ where $1 \leq h  \leq n^\delta$ and $p\in(0,1]$, all be real functions of $n$ which satisfy $aph \leq   \alpha\cdot \ln n $ and $bp >n^{\beta}$.  Then
		there exists  $C:=C(\alpha, \delta)\geq 1$ such that in $G(A,B,p)$, where $|A|=a $ and  $|B|=b$, we have	
			\begin{multline*}
				\Pr{\exists A'\subseteq A \text{ with }|A'|\leq \frac{1}{ph}\text{ and a dominating set $D\subset B$ of }A\setminus A'\text{ with }|D|\leq C\cdot \frac{\ln(eh)}{ph} }\\
				\geq  1- \exp\left[-\frac{n^{-\delta}\cdot b}{C}\right]. 
			\end{multline*}
	\end{restatable}

We believe that this result might be of independent interest, so we state and prove it separately in Section~\ref{sc:bip_dom}. The probability bound in this result allows us to establish necessary properties of neighbourhoods and their dominating sets, this is done in Section~\ref{sc:neighb_dom}. The proof of Theorem~\ref{th:fun_upper2} is completed in Section~\ref{sc:TH_upper_inter_proof}.

	\subsection{Neighbourhood Domination}
	\label{sc:neighb_dom}

In this section we apply our bound on the size of dominating sets (\cref{lem:small_t_dom_exist}) to find witnesses for functionality which are based on dominating the neighbours of a vertex. However, using just a single dominating set here is not enough. For a natural number $t$, a set $D \subseteq V \setminus S$ is \emph{$t$-dominating set} for $S$ if every vertex in $S$ has at least $t$ neighbours in $D$. Thus, a  $1$-dominating set is just a dominating set. We begin with \cref{lem:induced-domination}, which shows that for $t\geq 5$, \whp~$t$-dominating sets can be used to identify neighbourhoods in a random graph when $p$ is close to $p^*$.

All results in this section hold for $G_n\sim G(n,p)$ with 
\begin{equation}
\sqrt{\frac{\ln n}{n^{1+\varepsilon}}} \leq p \leq \sqrt{\frac{\ln n}{n^{1-\varepsilon}}}
\label{eq:bound_p_for_Sec_4_2}
\end{equation}
for some small $\varepsilon>0$.

\begin{claim}\label{lem:induced-domination}
	Let  $d\geq t\geq 5$. Then,  \whp~for every $y \in [n]$, every $S \subseteq N(y)$, and every $t$-dominating set $D$ for $S$ of size at most $d$, the number of vertices $w\notin N(y)$ such that, for some $z\in S$, $N_D(w)=N_D(z)$, is at most $4d$.  
	
\end{claim}
 
\begin{proof}	
	We will prove that  \whp~the following stronger event $\mathcal{B}$ holds. This event $\mathcal{B}$ states that for every $y \in [n]$ and every set $D$ of size at most $d$, there are at most $4d$ vertices $w\notin N(y)$ such that, 
	\begin{center}
		for some $z\in N(y)$ that has at least $t$ neighbours in $D$, $N_D(w)=N_D(z)$.
	\end{center} 
	The event $\mathcal{B}$ contains an intersection over many events, this makes it hard to bound $\mathbb{P}(\neg \mathcal{B})$ using a union bound. To overcome this we use a classic trick to bound difficult events. This is the more general observation that if there exists an `affable' event $\mathcal{A}$ and a `cover' event $\mathcal{C}$, both of whose probability we can bound easily, then provided $(\neg\mathcal{B}) \cap \mathcal{A}  \subseteq \mathcal{C}$, we have 
	\begin{equation}\label{eq:classictrick}
		\Pr{\neg\mathcal{B}} =  \Pr{(\neg\mathcal{B}) \cap \mathcal{A}  } + \Pr{\neg \mathcal{A} } \leq \Pr{\mathcal{C}}+ \Pr{\neg \mathcal{A} }.
	\end{equation}
Let $\mathcal{A}$ be the event that every vertex has degree at most $2np$, then $\Pr{\neg \mathcal{A}}=o(1)$~by  \cref{cl:degrees_common_degrees}. We now assume that the vertex set $[n]$ has an arbitrary but fixed ordering. Let $\mathcal{C}$ be the event there exists  $y \in [n]$, a set $D$ of size at most $d$,  and at least $4d+1$ vertices $w\notin N(y)$, such that 
\begin{center}
	one of the first $2np$ neighbours $z\in N(y)$ has at least $t$ neighbours in $D$ and $N_D(w)=N_D(z)$.
\end{center}  

Observe that $(\neg \mathcal{B} )\cap \mathcal{A}\subseteq \mathcal{C}$ holds, as if $\mathcal{A}$ holds then the first $2np$ neighbours of $y$ is all of $N(y)$. 

We will now bound $\mathbb{P}(\mathcal{C})$ by the union bound. To begin, the number of choices for $y$ is $n$, and for $D$ is at most $d{n\choose \lfloor d\rfloor}$. Now,  as soon as $y$ and $D$ are fixed and the first $2np$ vertices  $F\subseteq N(y)$ are exposed, then letting $\xi$ be the number of vertices $w\notin N(y)$ such that, for some $z\in F\cap N(D)$, $N_D(w)=N_D(z)$, we get that $\xi$ is stochastically dominated by $\mathrm{Bin}(n,2np\cdot p^t)$. This domination follows since there are less than $n$ choices for $w$, at most $2np$ choices for $z \in F$, and, as each $y$ has at least $t$ neighbours in $D$, the probability $w$ has the same neighbourhood in $D$ as $y$ is at most $p^t$. Thus, by applying Lemma~\ref{lm:Chernoff_AS} with $\lambda = \frac{4d}{2n^2p^{t+1}}>0$,
\[\mathbb{P}(\mathcal{C}) \leq n d{n\choose \lfloor d\rfloor} \cdot \mathbb{P}\left(\mathrm{Bin}(n,2np^{t+1})>4d\right)  \leq  
n^{d+2}\cdot  \exp\left[4d\ln\frac{2en^2p^{t+1}}{4d}\right]. \]Now, as $p \leq n^{-1/2 + \varepsilon}$ and $d\geq t\geq 5$, we have $\frac{2en^2p^{t+1}}{4d} \leq n^2p^{6} \leq  n^{-1 +6\varepsilon}\leq n^{-1/2}$. Thus, \[ \mathbb{P}(\mathcal{C}) \leq \exp\left[(d+2)\ln n-2d \ln n\right]=o(1).		\]
It now follows from \eqref{eq:classictrick} that $\mathbb{P}(\neg \mathcal{B})=o(1)$, as claimed. 
\end{proof}

For $\sqrt{\frac{\ln n}{n^{1+\varepsilon}}} \leq p \leq \sqrt{\frac{\ln n}{n^{1-\varepsilon}}}$, with some small $\varepsilon>0$, and a constant $C>0$, we define the function  
			\begin{equation}\label{eq:d}
d:=d(p,C)=  \left\{\begin{array}{cc}C\cdot \frac{np\ln(ew)}{\ln n}=C\cdot \frac{\ln(ew)}{pw}, & p=\sqrt{\frac{\ln n}{wn}},\, 1\leq w\leq n^{\varepsilon}; \\ C\cdot \frac{\ln(ew)}{p}, & p=\sqrt{\frac{w\ln n}{n}},\, 1\leq w\leq n^{\varepsilon} . \end{array}\right.\end{equation}

We can now adapt \cref{lem:small_t_dom_exist} to the setting of our upper bound. 
	\begin{claim} \label{cl:from_AB_reduction}
	 There exists $C>0$ such that    for $d:=d(p,C)$ given by \cref{eq:d} and any set $S\subset[n]\setminus\{1\}$ where  $|S|\leq 3np/2$, we have \begin{multline*}
			\mathbb{P}\left(\left.\exists D\subset[n]\setminus N[1]\,\,\exists N'\subset N(1):\,\,D\text{ is dominating for } N(1)\setminus N';\,|D|,|N'|\leq d\,\right|\, N(1)=S\right)
			\\
			\geq 1-\exp[-n^{1-\varepsilon}/C].
		\end{multline*}
 
	\end{claim}

	\begin{proof}
		Let us expose the neighbourhood of 1 and assume that $|N(1)|\leq	 3np/2$. Clearly, we may also assume that $|N(1)|\geq\ln^2 n$ --- otherwise, we can extend the set $N(1)$ arbitrarily and then a dominating set for the extension dominates $N(1)$ as well. The case $p=\sqrt{\frac{\ln n}{wn}}$ follows   from Theorem~\ref{lem:small_t_dom_exist} by taking $A=N(1)$, $B=[n]\setminus N[1]$, and $h=w$. To check the assumptions are satisfied,  we have $\ln^2 n \leq a:=|A|\leq 3np/2\leq b$, and also 
		\[
			ap\leq\frac{3}{2}np^2
			=  \frac{3\ln n}{2w} 
			= \frac{3}{2}\cdot\frac{\ln n}{h} \qquad \text{and} \qquad bp 
			\geq \Big(n-\frac{3np}{2}\Big)p 
			\geq n^{1/2-\varepsilon}.
		\]
		
		So, from now on, assume $p=\sqrt{\frac{w\ln n}{n}}$, where $1\leq w\leq n^{\varepsilon}$. Let 
		$$
		k=\frac{\ln(ew)}{p}.
		$$
		We divide $[n/2]\setminus N[1]$ into roughly equal parts $W_1,\ldots,W_m$, each of size either $\lfloor k\rfloor $ or $\lceil k \rceil$, where
		$$
		m=(1-o(1))\frac{n}{2k}.
		$$ 
 For each $j\in[m]$, let $\xi_j$ be the number of vertices in $N(1)$ that do not have neighbours in $W_j$. Then $\xi_j$ is a binomial random variable with 
		$$
		\mathbb{E}\xi_j\leq \frac{3}{2}\cdot np(1-p)^{\lfloor k\rfloor} =\Big(\frac{3}{2}  +o(1)\Big)\frac{np}{ew} \leq \frac{4}{5} \cdot \frac{np}{w} .
		$$
	We now apply Lemma~\ref{lm:Chernoff_JLR} with $t:= \frac{1}{5}\cdot \frac{np}{w}$ to bound 	$\mathbb{P}\left(\xi_j>\frac{np}{w}\right)$, whereby stochastic domination we can assume that  the estimate on $\mathbb{E}\xi_j$ above is an equality. Thus,
		$$
		\mathbb{P}\left(\xi_j>\frac{np}{w}\right)\leq \exp\left[-\frac{t^2}{2(	\mathbb{E}\xi_j +t/3)}\right] =\exp\left[-\frac{1/5^2}{ 2 (\frac{4}{5}+ \frac{1}{5\cdot 3} ) } \cdot \frac{np}{w}\right] \leq\exp\left[-\frac{np}{50w}\right].
		$$
		Since all $\xi_j$ are independent, and $w\leq n^\varepsilon$ the probability that each $W_j$ does not dominate more than $\frac{np}{w}$ vertices of $N(1)$ is at most 
		\begin{align*}
			\exp\left[-\frac{np}{50w}\cdot m\right] &= \exp\left[-(1-o(1))\frac{n^2p}{100wk}\right]\\
			&=\exp\left[-\frac{(1-o(1))n^2p^2}{100w\ln(ew)}\right]\\
			&=\exp\left[-\frac{(1-o(1))n\ln n}{100\ln(ew)}\right]\leq\exp\left[-\frac{n}{100}\right].
		\end{align*}
		Now, fix $j\in[m]$ such that $W_j$ dominates all but at most $\frac{np}{w}$ vertices of $N(1)$ and let $A$ be the set of vertices of $N(1)$ that are not dominated by $W_j$. If  $|A|< \ln^2 n$ then add some arbitrary vertices from $N(1)$ until $|A|\geq \ln^2 n$. We can now apply Theorem~\ref{lem:small_t_dom_exist} to this $A$ with $B= [n]\setminus([n/2]\cup N[1])$ and $h=1$, since $aph \leq  \frac{np}{w}\cdot p\cdot 1 = \ln n $ and $bp \geq n^{1/2}$, and we can take $\delta=\varepsilon/2>0$. This shows that there exists $C>0$ such that with probability at least $1-\exp[-n^{-\delta}b/C]\geq 1-\exp[- (1/2-o(1))n^{1-\varepsilon/2}/C]$, there exist $D\subset [n]\setminus([n/2]\cup N[1])$ of size at most $C/p$ and $N'\subset A$ of size at most $1/p$ such that $D$ dominates $A\setminus N'$. The set $D\sqcup W_j$ is the desired set (called $D$ in the statement), which dominates $N(1)\setminus N'$ where $|N'|\leq 1/p\leq d$, and this holds with the desired probability by the union bound.
	\end{proof}

	Let $\mathcal{Q}$ be the event that there are at least $n^{0.2}$ vertices among the first $\lceil n^{0.6}\rceil$ vertices of $G_n$  that have degrees at most $3np/2$. For a positive number $x$, let $\mathcal{Q}'(x)$ be the event that any two of the first $\lceil n^{0.6}\rceil$ vertices in $G_n$ have at most $x$ neighbours in common.

\begin{figure}[h]
	\begin{center}
		\includegraphics[scale=1]{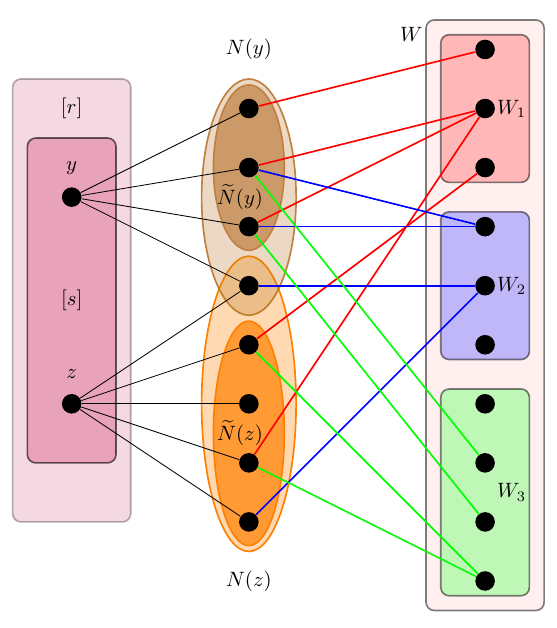}
	\end{center}
\caption{Illustration of the proof of \cref{lem:clm1}. The vertex $y\in [s]$ has a $3$-dominating set for two of the vertices in its neighbourhood.   }\label{fig:boostclmfig}
\end{figure}
\begin{claim} \label{lem:clm1}
Let $x=o(n^{0.7}\cdot p)$. Then, there exists a $C>0$ such that for $d:=d(p,C)$ given by \eqref{eq:d}, any $t\in\mathbb{Z}_{>0}$, and   
				\begin{equation}
					\mathcal{A}=\left\{\exists y\,\, \exists D\subset[n]\setminus N[y]\,\,\exists N'\subset N(y):\,\,D\text{ }t\text{-dominates } N(y) \setminus N';\,|D|,|N'|\leq 2td\right\},
					\label{eq:A_def}
				\end{equation}
we have $\mathbb{P}\left(( \neg \mathcal{A}) \cap \mathcal{Q}\cap \mathcal{Q}'(x)\right) \leq\exp[-n^{1.1}].$
\end{claim}

			\begin{proof}
Let $\mathcal{Q}':=\mathcal{Q}'(x)$. For the event $\mathcal{A}$ given by \eqref{eq:A_def} we have
				$$
				\mathbb{P}((\neg \mathcal{A})\cap \mathcal{Q}\cap \mathcal{Q}')\leq
				\mathbb{P}( \neg \mathcal{A}\mid \mathcal{Q}\cap \mathcal{Q}').
				$$ 
Thus, it suffices to show that $\mathbb{P}(\neg \mathcal{A}\mid \mathcal{Q}\cap \mathcal{Q}')  \leq \exp[-n^{1.1}]$, which we will do.

Let $r=\lceil n^{0.6}\rceil$ and $s=\lceil n^{0.2}\rceil$. Expose the neighbours of vertices $1,\ldots,r$. If $\mathcal{Q}$ holds, then assume without loss of generality that the sets $N(1),\ldots,N(s)$ each have size at most $3np/2$. Let $W=[n]\setminus (N(1)\cup\ldots\cup N(s)\cup[r])$. Partition $W$ into $t$ parts $W_1,\ldots,W_t$ with sizes as equal as possible. Then each $W_\ell$, for $\ell \in [t]$, has $\frac{n}{t}(1-o(1))$ vertices, see \cref{fig:boostclmfig}. 

For $i\in[s]$, we define $\widetilde{N}(i)$ to be the part of the $i$'s neighbourhood not contained $N(j)$ for any other $j\neq i\in [s]$, that is 
\[
\widetilde{N}(i) := N(i)\setminus \bigcup_{j\in [s], j\neq i} N(j). 
\]
Observe that conditional on $\mathcal{Q}\cap \mathcal{Q}'$ we have $|\widetilde{N}(i)|=(1+o(1))|N(i)|$ for all $i\in [s]$. Additionally, for  $i\in[s]$, let $\mathcal{E}_i$ be the event that all the sets $W_1,\ldots,W_t$ contain a dominating set $D$ of $\widetilde{N}(i)\setminus N'$ for some $N'\subset \widetilde{N}(i)$ where $|D|, |N'|\leq d$. Take $S_1,\ldots,S_s\subset[n]$ to be any sets such that the event
$$
 \mathcal{N}_{S_1,\ldots,S_s}:=\{N(1)=S_1,\ldots,N(s)=S_s\}
$$ 
implies $\mathcal{Q}\cap \mathcal{Q}'$.  To bound $\mathbb{P}\left(\neg \mathcal{E}_i\mid \mathcal{N}_{S_1,\ldots,S_s}\right)$ we will first bound the probability that, for a given set $W_\ell$, where $\ell\in[t]$, there exist the desired sets $D, N'$ such that $D\subset W_\ell$ dominates $\widetilde{N}(i)\setminus N'$. If we take $C>0$ suitably large, then it follows from applying \cref{cl:from_AB_reduction} to the random graph on $\{i\}\cup \widetilde{N}(i)\cup W_\ell$ that the probability no such $D,N'$ exist is at most $\exp\Big[-\frac{  ((1-o(1))n/t)^{1-\varepsilon}}{C}\Big]$. Thus by the union bound, we have that for any $i\in[k]$
$$
\mathbb{P}\left(\neg \mathcal{E}_i\mid \mathcal{N}_{S_1,\ldots,S_s}\right)\leq t\cdot \exp\Big[-\frac{  ((1-o(1))n/t)^{1-\varepsilon}}{C}\Big]\leq\exp[-n^{0.9}].
$$

Observe that if  $\mathcal{Q}'$ holds then $|N(i)\setminus \widetilde{N}(i)| \leq s \cdot x =o(n^{0.2}\cdot n^{0.7} p)=o(d)$ for all $i\in [s]$, as $\varepsilon >0$ is small. Thus, if additionally $\mathcal{E}_i$ holds, for some $i\in[s]$, then, by adding up all the dominating sets from each $W_\ell$, we have a $t$-dominating set of $N(i)$ of size $t|D|\leq td$, and we can add $|N(i)\setminus \widetilde{N}(i)|$ along with all the sets $N'$ to give a subset $N'' \subseteq N(i)$ (this is $N'$ in the statement) of vertices that are not covered with $|N''|\leq  t|N'| + |N(i)\setminus \widetilde{N}(i)|\leq 2td$. Note  that for any $i,j\in [s]$, conditional on the event $\mathcal{N}_{S_1,\ldots,S_s}$, the events $\neg \mathcal{E}_i$ and $\neg \mathcal{E}_j$ are independent as $\widetilde{N}(i)\cap \widetilde{N}(j)=\emptyset$, i.e.~the overlapping parts of the neighbourhoods of $i$ and $j$ are not considered within these events. Thus, 
\begin{align*}
\mathbb{P}(\neg \mathcal{A}\mid \mathcal{Q}\cap \mathcal{Q}')&\leq \mathbb{P}(\neg \mathcal{E}_1\wedge\ldots\wedge\neg{\mathcal{E}_{s}}\mid \mathcal{Q}\cap \mathcal{Q}') \\
&=
\sum_{S_1,\ldots,S_s}\frac{\mathbb{P}(\neg{\mathcal{E}_1}\wedge\ldots\wedge\neg{\mathcal{E}_{s}}\mid \mathcal{N}_{S_1,\ldots,S_s})\mathbb{P}(\mathcal{N}_{S_1,\ldots,S_s})}{\mathbb{P}(\mathcal{Q}\cap \mathcal{Q}')}\\
&= \sum_{S_1,\ldots,S_s}\frac{\mathbb{P}(\neg{\mathcal{E}_1}\mid \mathcal{N}_{S_1,\ldots,S_s})\ldots\mathbb{P}(\neg{\mathcal{E}_{s}}\mid \mathcal{N}_{S_1,\ldots,S_s})\mathbb{P}(\mathcal{N}_{S_1,\ldots,S_s})}{\mathbb{P}(\mathcal{Q}\cap \mathcal{Q}')}\leq\exp[-s\cdot n^{0.9}].
						\end{align*}
						This completes the proof.
					\end{proof}

\subsection{Proof of Theorem~\ref{th:fun_upper2}}
			\label{sc:TH_upper_inter_proof}
			
Let $\varepsilon>0$ be sufficiently small so that for all $p$ satisfying~\eqref{eq:bound_p_for_Sec_4_2}, all statements from Section~\ref{sc:neighb_dom} hold true. 
 We can assume w.l.o.g.\ that $\varepsilon <0.001$, and then fix the epsilon in the statement of \cref{th:fun_upper2} to be $\varepsilon/8$. Thus in terms of our $\varepsilon>0$, we have \begin{equation} \label{eq:peps}
	p\in [n^{-1/2-\varepsilon/8},n^{-1/2+\varepsilon/8}].
\end{equation} 		
Now, let $t=5$ and $C>0$ be a sufficiently large constant, and for this $C$ let $d:=d(p,C)$ be as given by \eqref{eq:d}. To bound the functionality of $G_n \sim G(n,p)$ from above we must bound the functionality in all induced subgraphs. We break into two cases depending on $k$, the size of the induced subgraph. 

\paragraph{Large ($k\geq n^{1-\epsilon/3}$).} We denote by $\mathcal{A}$ the event from Claim~\ref {lem:clm1}, given by~\eqref{eq:A_def}. For $U\subset[n]$, let $\mathcal{Q}_U$ be the event that there are at least $|U|^{0.2}$ vertices among the first $\lceil |U|^{0.6}\rceil$ vertices of $G_n[U]$  that have degrees at most $3|U|p/2$, i.e.~$\mathcal{Q}_U=\big\{G_n[U]\in \mathcal{Q}\big\}$, recalling that $\mathcal{Q}$ and $\mathcal{Q}'(x)$ were defined in Section~\ref{sc:neighb_dom}.  In the same way, we set $\mathcal{Q}'_U=\big\{G_n[U]\in \mathcal{Q}'(k^{0.65}\cdot p)\big\}$ and  $\mathcal{A}_U=\{G_n[U]\in \mathcal{A}\}$.  We also let $\mathcal{Q}''$ be the event that any two vertices in $G_n$ have at most $n^{0.6}\cdot p$ neighbours in common. Notice that $\mathcal{Q}''$ implies $\mathcal{Q}'_U$ for any $U$ with $|U|\geq n^{1-\varepsilon/3}$. Finally, since $n^{0.6}p\geq \max\{2\ln n,10np^2\}$, Lemma~\ref{cl:degrees_common_degrees} implies that 
\begin{equation}\label{eq:PnegQ} \mathbb{P}(\neg \mathcal{Q}'')=o(1).
\end{equation}

Our aim in this case is to show that \whp for all $U$ with $|U|\geq n^{1-\varepsilon/3}$, in the induced subgraph $G_n[U]$ there is $y\in U$, $D\subset U\setminus N[y]$ and some $N'\subset N(y)$ such that $D$ $t$-dominates $N(y) \setminus N'$ and $|D|,|N'|\leq 2td$, where $d=d(p(n),C)$. That is, \whp for all such $U$ the event $G_n[U] \in \mathcal{A} $ holds. We will achieve this by applying a union bound over such sets $U$, then use \cref{lem:clm1} to bound the probability $G_n[U] \in \mathcal{A}$. There is a small issue to overcome here as each $ G_n[U]\sim G(|U|, p(n))$, however \cref{lem:clm1} must be applied to a random graph from $G(n,p(n))$. Let $G_k\sim G(k, p(n))$, where $n^{1-\varepsilon/3}\leq k\leq n$. Note that w.l.o.g. $k$ is a bijective function of $n$, so there is a function $p'(k)\in[0,1]$ such that, for each $n \in \mathbb{N}$ and $k=k(n)$, $ p'(k) =p(n)$ and thus for this $p'$ we have $G_k\sim G(k,p'(k))$. As $ n^{1-\varepsilon/3}\leq k \leq n$ corresponds to $k\leq n \leq k^{\frac{1}{1-\varepsilon/3}} $ and we can assume $\varepsilon < 0.001$, we have 
\[
 \sqrt{\frac{\ln k}{k^{1+\varepsilon}}} \leq k^{-\frac{1}{2}- \frac{5\varepsilon}{12-4\varepsilon }}  = \left(k^{\frac{1}{1-\varepsilon/3}}\right)^{-1/2-\varepsilon/4} \leq     p'(k) =p(n)\leq   k^{-1/2+\varepsilon/4} \leq \sqrt{\frac{\ln k}{k^{1-\varepsilon}}}.     
 \] 
 Since $\varepsilon >0$ is small enough, we can apply \cref{lem:clm1} to $G_k\sim G(k,p'(k))$ and obtain \begin{equation}\label{eq:prg_k} \mathbb{P}(G_k\in \neg \mathcal{A}\cap \mathcal{Q}\cap \mathcal{Q}') \leq \exp[-k^{1.1}].\end{equation} 
The bound has a good rate of decay, however the event $G_k\in \mathcal{A}$ gives  $|D|,|N'|\leq 2td(p'(k),C)$, whereas we want $|D|,|N'|\leq 2td(p(n),C)$. The following claim shows this is not an issue.

\begin{claim}Let $C>0$ and $d:=d(p,C)$ be given by \eqref{eq:d}. Then for any $p=p(n)\in[0,1]$, all $n$ large enough, and $n^{1-\varepsilon/3} \leq k\leq n$,  \[d(p'(k),C)\leq d(p(n),C).\]  
\end{claim}
\begin{proof}
	First suppose that $ p(n)=\sqrt{\frac{\ln n}{wn}}=p'(k)$ for some $1\leq w\leq n^{\varepsilon}$. Then, since we can assume $k,n$ are large,  $\frac{\ln n}{n}\leq \frac{\ln k}{k}$ and thus there is a $w'\geq w\geq 1$ such that $ p'(k)=\sqrt{\frac{\ln k}{w'k}}$. Now since $\frac{\ln (ex)}{x}$ is decreasing in $x\geq 1$, in this case,  
	\[
	d(p'(k),C) = \frac{C\ln(ew')}{p'(k)w'} \leq \frac{C\ln(ew)}{p(n)w}=  d(p(n),C) .
	\]    
	Otherwise, if $ p(n)=\sqrt{\frac{w\ln n}{n}}=p'(k)$ for some $1\leq w\leq n^{\varepsilon}$, then there exists a $n^{-\varepsilon}\leq w'\leq w$ such that $ p'(k)=\sqrt{\frac{w'\ln k}{k}}$. First suppose that $1\leq w'\leq w$, then 
	$$
	 d(p'(k),C) = \frac{C\ln(ew')}{p'(k)} \leq \frac{C\ln(ew)}{p(n)}= d(p(n),C).
	$$
	Otherwise,  $n^{-\varepsilon}\leq w'<1 $, and then $ p'(k)=\sqrt{\frac{\ln n}{(1/w') n}}$, where $1/w'>1$. Therefore,
	 \[
	 d(p'(k),C) = \frac{C\ln(e/w')}{p'(k)/w'} \overset{(\star)}{\leq} \frac{C}{p'(k)} \overset{(\dagger)}{\leq}   \frac{C\ln(ew)}{p(n)}= d(p(n),C),
	 \] 
	 where $(\star)$ holds since $\frac{\ln(ex)}{x} $ is decreasing in $x\geq 1$, and $(\dagger)$ holds since $w\geq 1 $.     
\end{proof}

  Now, by \cref{eq:PnegQ} and the above discussion around  $G_k \in \mathcal{A}$ vs $G_n[U] \in \mathcal{A}_U$, we have
\begin{align}
	\mathbb{P}(\exists U\subset[n]:&\, |U|\geq n^{1-\varepsilon/3}\text{ and }\neg \mathcal{A}_U) \notag \\
	 &\leq \mathbb{P}(\neg \mathcal{Q}'')+\mathbb{P}(\exists U\subset[n]:\, |U|\geq n^{1-\varepsilon/3}\text{ and }\neg \mathcal{A}_U\cap \mathcal{Q}'')\notag \\
	&\leq o(1)+\mathbb{P}(\exists U\subset[n]:\, |U|\geq n^{1-\varepsilon/3}\text{ and }\neg \mathcal{A}_U\cap \mathcal{Q}'_U)\notag \\
	&\leq o(1)+\sum_{k\geq n^{1-\varepsilon/3}}\binom{n}{k}\biggl(\mathbb{P}(G_k\in \neg \mathcal{A}\cap \mathcal{Q}\cap \mathcal{Q}')+\mathbb{P}(G_k\notin\mathcal{Q})\biggr).\label{eq:4.1first}
\end{align}   
We are now in a position to apply \cref{eq:prg_k}, but we will also need to following to bound the other probability on the RHS of \cref{eq:4.1first}. 
\begin{claim}\label{clm:SecondRhs}
	For every $k\geq n^{1-\varepsilon/3}$, $\mathbb{P}(G_k\notin\mathcal{Q})=o(2^{-n})$.
\end{claim}			
\begin{proof}
Recall that we assume $\varepsilon<0.001$. Let $S$ be the set of the first $\lceil k^{0.6}\rceil$ vertices in $[k]$. By Lemma~\ref{lm:Chernoff_JLR}, with probability $1-o(2^{-n})$ there are $o(k^{1.6}p)$ edges in $G_k[S]$. Similarly, with probability $1-o(2^{-n})$ there are at most $1.1 k^{1.6}p$ edges between $S$ and $[k]\setminus S$. However, the event $G_k\notin\mathcal{Q}$ implies that almost all vertices in $U$ have degrees at least $\frac{3}{2}kp$ in $G_k$, which is only possible when either $G_k[U]$ has more than $0.1 k^{1.6}p$ edges or there are more than $1.1k^{1.6}p$ edges between $S$ and $[k]\setminus S$. This completes the proof.
\end{proof}

Now, inserting \cref{eq:prg_k} and \cref{clm:SecondRhs} into \cref{eq:4.1first} gives
 \begin{align*} 
	\mathbb{P}(\exists U\subset[n]:\, |U|\geq n^{1-\varepsilon/3}\text{ and }\neg \mathcal{A}_U)	&\leq o(1)+\sum_{k\geq n^{0.99}}{n\choose k}\exp[-k^{1.1}]\\
	&\leq o(1)+\sum_{k\geq n^{0.99}}2^n \exp\left[-n^{0.99\cdot 1.1}\right]=o(1).
\end{align*}
 
			So, \whp~in any induced subgraph $H$ of size at least $n^{1-\varepsilon/3}$, there exists a vertex $y$, a set $N'\subset N(y)\cap V(H)$ of size at most $2td$, and a $t$-dominating set of $N(y)\cap V(H)\setminus N'$ of size at most $2td$. Let us fix such a set $D:=D(H;y)\subset V(H)$, which $t$-dominates the set $S:=N(y)\cap V(H)\setminus N'$. By Claim~\ref{lem:induced-domination}, we can assume that the number of vertices $w$ in $H$ that do not belong to $N(y)$ such that, for some $z\in S$, $N_D(w)=N_D(z)$, is at most $4\cdot 2td$. Denote the set of these vertices by $B$ and note that, by \cref{obs:fun1}, $y$ is a function of $D \cup N' \cup B$. Thus, \whp
			$$
			\max_{H\sqsubset G_n:\,V(H)\geq n^{1-\varepsilon/3}}\;\min_{y\in V(H)}\fun(y)\leq 2td +2td + 8td = 60d.
			$$
		\paragraph{Small ($k\leq n^{1-\epsilon/3}$).}	 
		For this case we will prove the stronger bound: 	\whp, 
		$$
		\max_{H\sqsubset G_n:\,V(H)< n^{1-\epsilon/3}}\;\min_{y\in V(H)}\fun(y)\leq n^{1/2-\varepsilon/7} <d.
		$$
		This follows immediately from the next claim, completing the proof of Theorem~\ref{th:fun_upper2}.

				\begin{claim}
				Let $\varepsilon>0$ be sufficiently small,  $n^{-1/2-\varepsilon/8}\leq p\leq n^{-1/2+\varepsilon/8}$, and $G_n\sim G(n,p)$. Then, \whp, every $H\sqsubset G_n$ with $V(H)< n^{1-\varepsilon/3}$ has degeneracy at most $n^{1/2-\varepsilon/7}$.
			\end{claim}
			
			\begin{proof}
				Let $k<n^{1-\varepsilon/3}$. Then, by Lemma~\ref{lm:Chernoff_AS} with $\lambda=\frac{n^{1/2-\varepsilon/7}k}{k^2p}$, we have
				\begin{equation*}
					\mathbb{P}\left(|E(G_k)|>n^{1/2-\varepsilon/7}k\right)\leq
					\mathbb{P}\left(\mathrm{Bin}(k^2,p)>n^{1/2-\varepsilon/7}k\right)\leq
					\exp\left[n^{1/2-\varepsilon/3}k\cdot \ln\frac{e k^2p}{n^{1/2-\varepsilon/7}k}\right].
				\end{equation*}
Then, as 
$$
\ln\frac{ekp}{n^{1/2-\varepsilon/7}} \leq  1+\left(1-\frac{\varepsilon}{3} -\left(\frac{1}{2}-\frac{\varepsilon}{8}\right) -\left(\frac{1}{2}-\frac{\varepsilon}{7}\right)\right)\ln  n +1  = 1-\frac{11}{168}\varepsilon\cdot\ln  n,
$$
for large $n$ we have 
				\begin{equation*}
					\mathbb{P}\left(|E(G_k)|>n^{1/2-\varepsilon/7}k\right) \leq \exp\left[-\frac{\varepsilon}{20}\cdot n^{1/2-\varepsilon/7}k\ln n\right]. 
				\end{equation*} 
				By the union bound and as $\binom{n}{k} \leq \exp(k\ln n)$, the probability that there exist   $H\sqsubset G_n$ with $V(H)< n^{1-\varepsilon/3}$ and degeneracy at least $n^{1/2-\varepsilon/7}$ is at most
				\begin{align*}
					\sum_{k=1}^{n^{1-\varepsilon/3}}{n\choose k}\exp\left[-\frac{\varepsilon}{20}\cdot n^{1/2-\varepsilon/7}k \ln n\right]&\leq 
					\sum_{k=1}^{n^{1-\varepsilon/3}}\exp\left[k\ln n-\frac{\varepsilon}{20}\cdot n^{1/2-\varepsilon/7}k \ln n\right]\\
					&= n^{ -(\varepsilon/20+o(1))\cdot n^{1/2-\varepsilon/7} }=o(1),
				\end{align*}
				which completes the proof since functionality is at most degeneracy by \cref{obs:fun-dgn}.
			\end{proof}

			\section{Lower Bounds}\label{sec:lower-bounds}
			By \cref{obs:fun1}, a vertex $y \in V(G)$ is not a function of a set $S \subseteq V(G) \setminus \{y\}$ if there exists a pair of vertices $z_1,z_2 \notin S \cup \{y\}$ such that $z_1$ and $z_2$ have the same neighbourhood in $S$ but different adjacencies with $y$. 
			
			We apply this idea to prove our lower bounds by showing that such a pair $z_1,z_2$ exists for all such $y$, and sets $S$ of a certain size.  We split the proof into two cases for small and large $p$, Sections \ref{sec:lowersmallp} and \ref{sec:lowerlargep} respectively, as the way in which this idea is applied is dependent on $p$. 
			
			\subsection{Lower Bound for \texorpdfstring{$p\leq p^*$}{p at most p*}}\label{sec:lowersmallp}

			\begin{lemma}\label{lem:lower-p-less-p-star}
				Let $G_n \sim G(n,p)$, where $p=\sqrt{\frac{\ln n}{wn}}$, $w=w(n)\geq 1$.
				Then  \whp
				\[
				\fun\left(G_n\right) \geq  \left\lfloor\frac{np\ln(ew)}{20\ln n}\right\rfloor. 
				\] 
			\end{lemma}
			\begin{proof}
				Without loss of generality we can assume that $np\geq 19$, since otherwise $\frac{np\ln(ew)}{20\ln n}<1$, and so the statement holds trivially. This can be seen for  $w\leq n^{1.02}$ as $\ln(ew)$ becomes sufficiently small, and for $w>n^{1.02}$ as the numerator $np\ln(ew)$ becomes sufficiently small. Note that the assumption $np \geq 19$ implies that $w \leq n \ln n$.

				Consider the following process. Iteratively, remove from $G_n$ vertices of degree at most $np/5$. Let $k$ be the number of removed vertices, and $H_n$ be the remaining graph on $n_0 := n-k$ vertices. Let 
				\[\mathcal{D}:= \Big\{ \Delta \leq\max\{\ln n,4np\} \leq \frac{np\ln n}{19}\Big\},\] and define the following two events 
				\begin{equation}
					\mathcal{A}:=\big\{n_0 \geq n / \ln n\big\}\cap \mathcal{D}, \qquad \text{and}\qquad	\mathcal{A}_{\mathrm{large}}:=\big\{n_0=n\big\}\cap \mathcal{D},   \end{equation} and note that $\mathcal{A}_{\mathrm{large}}\subseteq \mathcal{A}$. We now bound the probability of these events.
				\begin{claim}\label{clm:lower1}
				$\mathbb{P}(\mathcal{A}) =1-o(1)$. Furthermore,	if $\frac{4\ln n}{n}\leq p\leq\frac{1}{2}$, then $\mathbb{P}(\mathcal{A}_{\mathrm{large}}) =1-o(1)$.
				\end{claim}
				
				\begin{pocd}{\ref{clm:lower1}}
				Due to the union bound and Chernoff bound, \whp $\mathcal{D}$ holds. 
				
				The second part of the claim is straightforward since all degrees in $G(n,p)$ are at least $np/5$ \whp~when $p\geq\frac{4\ln n}{n}$ by the union bound and Chernoff bound (Lemma~\ref{lm:Chernoff_JLR}). 
			
			For the first part, observe that the total number of edges in $G_n$ is at most $k n p/5 + n_0 \Delta$, where $\Delta$ is the maximum degree in $G_n$. 
					Thus, since \whp~$|E(G_n)| = \binom{n}{2}p(1 + o(1))$,
					we have that \whp
					\[
					\binom{n}{2}p(1 +o(1)) \leq \frac{knp}{5} + n_0\Delta \leq 
					\frac{1}{5}np\left( k + n_0\ln n \right).
					\]
					Thus, for $n$ large enough, we have $2n \leq k + n_0\ln{n}$, which together with the fact that $k \leq n$ implies the desired 	$n_0 \geq n/\ln n$. 
				\end{pocd}

				In the rest of the proof we will show that 
				\whp the functionality of $H_n$ is at least 
				$$
				m:=\left\lfloor\frac{np\ln(ew)}{20\ln n}\right\rfloor.
				$$
				By definition, this will imply the desired bound for the functionality of $G_n$.
				To do so we will use \cref{obs:fun1} to prove that for any set $S \subseteq V(H_n)$ with $|S| = m$ and any vertex $y$ outside $S$, the vertex $y$ is not a function of $S$. Specifically, we will show that for any such $S$ and $y$ there exist $z_1,z_2 \not\in S \cup \{y\}$ such that $z_1$ and $z_2$ have different adjacencies with $y$, and neither $z_1$ nor $z_2$ has neighbours in $S$.
				
				Fix $S\subset[n]$ of size $m$,  $y\in[n]$, and denote by $N$ the neighbourhood of $y$ in $H_n$ (though we are only interested in $S\subset V(H_n)$ and $y\in V(H_n)$, we want to treat $S$ and $y$ as deterministic objects, thus we choose them from the entire set of vertices).  Assuming that $y$ belongs to $H_n$ and $\mathcal{A}$ holds, we have $|V(H_n)|\geq n/\ln n$, and $\Delta\leq\frac{np\ln n}{19}$, and since $m \leq \left\lfloor (1+o(1))\frac{np}{20} \right\rfloor$ we get that
				\begin{equation} \label{eq:Z_sizes_bound1}
					|N\setminus S|\geq\frac{np}{5}-m\geq\frac{np}{7},
				\end{equation}
				and
				\begin{equation} \label{eq:Z_sizes_bound2}
					|V(H_n)\setminus(N\cup S\cup\{y\})|\geq\frac{n}{\ln n}-\frac{np\ln n}{19}-m-1=(1-o(1))\frac{n}{\ln n}.
				\end{equation}
				Letting 
				$Z\subseteq[n]\setminus S$ be any set, we have 
				\begin{align}
					\mathbb{P}(\text{every }z \in Z \text{ has a neighbour in } S)
					& =(1-(1-p)^{m})^{|Z|}\leq(mp)^{|Z|}.
					\label{eq:Z-S}
				\end{align}

				Let us now separately consider two cases: $\frac{19}{n}\leq p<\frac{4\ln n}{n}$ and $\frac{4\ln n}{n}\leq p\leq\sqrt{\frac{\ln n}{n}}$.
				
				\paragraph{Small $p$.} Let $\frac{19}{n}\leq p<\frac{4\ln n}{n}$. For this range of $p$, we have $w\in [n/(16\ln n),(n\ln n)/361]$, and thus  
				\begin{equation}
					m=\left\lfloor(1+o(1))\frac{np}{20}\right\rfloor\quad\text{ and }\quad mp\leq \frac{(\ln n)^2}{n}.
					\label{eq:small_m_bounds}
				\end{equation}
				We start by showing  $\mathbb{P}(\mathcal{B}_1)=1-o(1)$, where $\mathcal{B}_1$ is the event that for every $y\in[n]$, $S\subset[n]\setminus\{y\}$ of size $m$, and $N\subset N(y)$ of size at least $np/5$, there exists $z_1 \in N \setminus S$ with no neighbours in $S$. 
				
				We will use the trick \eqref{eq:classictrick} from \cref{lem:induced-domination}, and define an event $\mathcal{C}_1$ satisfying $(\neg\mathcal{B}_1) \cap \mathcal{A} \subseteq \mathcal{C}_1 $. It then suffices to show $\mathbb{P}(\mathcal{C}_1)=o(1)$, as $\mathbb{P}(\neg \mathcal{A})=o(1)$ by \cref{clm:lower1}. To define $\mathcal{C}_1$ we first fix an ordering of $[n]$. Now assign, for each vertex $y$, a set $F_y$ containing the first (at most) $\lceil 16 \ln n \rceil$ vertices of $N(y)$ and, if $|N(y)|<\lceil np/5 \rceil$, then $F_y$ also contains the first $\lceil np/5 \rceil -|N(y)|$ vertices in order that are not present in $N(y)$. Let $\mathcal{C}_1$ be the event that there exist $y\in[n]$, $S\subset[n]\setminus\{y\}$ of size $m$, and $N\subset F_y$ of size exactly $\lceil np/5 \rceil$, where every vertex in $N \setminus S$ has a neighbour in $S$. Now, by  \eqref{eq:Z_sizes_bound1}, \eqref{eq:Z-S}, and the union bound,
					\begin{align}
				\mathbb{P}(\mathcal{C}_1) \leq	n \cdot \binom{n}{m} \cdot \binom{ \lceil 16\ln n\rceil}{\lceil np/5\rceil} \cdot (mp)^{np/7}
					& \leq n^{m+1} \left(mp(\ln n)^3\right)^{np/7}.
					\label{eq:small_union_bound}
				\end{align}
				Since $np\geq 19$,  due to~\eqref{eq:small_m_bounds}, for large $n$ we have that     
				$$
				n^{m+1} \left(mp(\ln n)^3\right)^{np/7}
				\leq 
				n^{m+1} \left(\frac{(\ln n)^2}{n}\cdot (\ln n)^3\right)^{19/14} \cdot \left(mp(\ln n)^3\right)^{np/14}
				\leq  n^{m} \left(mp(\ln n)^3\right)^{np/14}.
				$$
				Thus, again, as $n^m \leq \exp[\frac{1}{19}np\ln n] $  due to~\eqref{eq:small_m_bounds}, the desired probability satisfies
				\begin{align}
						\mathbb{P}(\mathcal{C}_1)  
						& \leq\exp\left[\frac{1}{19}np\ln n+\frac{1}{14}np\ln\left(mp(\ln n)^3\right)\right] =
					\exp\left[\frac{1}{19}np\ln n-\frac{1-o(1)}{14}np\ln n\right]=o(1).\label{eq:boundon6}
				\end{align}

				Next we show  $\mathbb{P}(\mathcal{B}_2)=1-o(1)$, where $\mathcal{B}_2$ is the event that for every set $U\subset[n]$ of size at least $n/\ln n$, $y\in U$, $S\subset U\setminus\{y\}$ of size $m$, and $N\subset N(y)\cap U$ of size at least $np/5$ there exists  
				$z_2 \in U\setminus(N\cup S\cup\{y\})$ with no neighbours in $S$. We show this by setting up a similar event $\mathcal{C}_2$ using sets $F_y$. The event $\mathcal{C}_2$ states that there exist $U\subset[n]$ of size at least $n/\ln n$, $y\in U$, $S\subset U\setminus\{y\}$ of size $m$, and $N\subset F_y\cap U$ of size exactly $\lceil np/5 \rceil$, such that every vertex in $U\setminus(N\cup S\cup\{y\})$ has a neighbour in $S$.
				Indeed, from \eqref{eq:Z_sizes_bound2} and \eqref{eq:Z-S}, by the union bound,  
			
				\begin{equation*}
					\mathbb{P}(\mathcal{C}_2) \leq 2^n\cdot n \cdot \binom{n}{m} \cdot \binom{ \lceil 16\ln n\rceil }{ \lceil np/5\rceil} \cdot (mp)^{(1-o(1))n/\ln n}. \end{equation*}
				Observe that bound on 	$\mathbb{P}(\mathcal{C}_2)$ contains the first bound~\eqref{eq:small_union_bound} we obtained on $\mathbb{P}(\mathcal{C}_1)$ as a factor. Thus applying the bound \eqref{eq:boundon6} that we have just obtained for~\eqref{eq:small_union_bound} to this factor  yields
				\begin{equation*}  	\mathbb{P}(\mathcal{C}_2)   \leq 2^n\cdot o(1)\cdot  \left(mp\right)^{(1-o(1))n/\ln n - np/7}
				 \stackrel{\eqref{eq:small_m_bounds}}\leq 2^n \left(  (\ln n)^2/n\right)^{(1-o(1))n/\ln n} = 2^n e^{-(1-o(1))n}=o(1).
				\end{equation*}

				Finally, let $\mathcal{B}$ be the event that for every choice of $y\in V(H_n)$ and $S\subset V(H_n)\setminus y$ of size $m$, there exist $z_1 \in N(y)\cap V(H_n) \setminus S$ and $z_2 \in V(H_n) \setminus (N(y) \cup S \cup \{y\})$ that both do not have neighbours in $S$. By \cref{obs:fun1} and the definition of functionality, if $\mathcal{B}$ holds then $\fun(G_n) \geq \fun(H_n) \geq m$ as desired. However, $\mathcal{B} \supseteq\mathcal{B}_1\cap \mathcal{B}_2 $ and thus $\mathbb{P}(\neg \mathcal{B}) \leq 	\mathbb{P}(\mathcal{C}_1) + 	\mathbb{P}(\mathcal{C}_2) + 	\mathbb{P}(\neg\mathcal{A}) =o(1)  $.

				\paragraph{Large $p$.} Let $\frac{4\ln n}{n}\leq p\leq\sqrt{\frac{\ln n}{n}}$. In this case $w\in [1,n/(16\ln n)]$, and thus as $p=\sqrt{\frac{\ln n}{wn}}$, we have
				\begin{equation}
					m=(1+o(1))\frac{np\ln(ew)}{20\ln n}\quad\text{ and }\quad mp\leq(1+o(1))\frac{\ln(ew)}{20w}.
					\label{eq:large_m_bounds}
				\end{equation}
		Similarly to before, we give an arbitrary order to the vertices of $G_n$ and assign for each vertex $y$, a set $F_y$ such that either $F_y=N(y)$ if $|N(y)|\geq\lceil np/5 \rceil$ or, otherwise, $F_y$ contains all vertices from $N(y)$ and also the first $\lceil np/5 \rceil -|N(y)|$ vertices in order that are not present in $N(y)$. Let $\mathcal{C}$ be the event that there exist $y\in[n]$ and $S\subset[n]\setminus\{y\}$ of size $m$ such that either every $z_1 \in F_y \setminus S$ has a neighbour in $S$, or every $z_2\in[n]\setminus(F_y \cup S\cup\{y\})$ has a neighbour in $S$. We now bound  $\mathbb{P}(\mathcal{C})$ by the union bound, using the fact that both $F_y$ and $[n]\setminus(F_y \cup S\cup\{y\})$ have size at least $np/5$ while $m<np/19$,  which gives
			\[\mathbb{P}(\mathcal{C}) \leq 	2\cdot n \cdot \binom{n}{m} \cdot (mp)^{np/7}  \leq n^{m+1} \left(mp\right)^{np/7}. \] 
			Now, by \eqref{eq:large_m_bounds} and using the inequality $(*)$ that $20w>(\ln(ew))^{2}$ for all $w\geq 1$, we have 
				\begin{align}
					\mathbb{P}(\mathcal{C}) &\stackrel{\eqref{eq:large_m_bounds}}\leq\exp\left[(1+o(1))\frac{np\ln(ew)}{20}-(1+o(1))\frac{np\ln(20w/\ln(ew))}{7}\right]\notag \\
					&\stackrel{(*)}\leq\exp\left[\left(\frac{1}{140}+o(1)\right)(7np\ln(ew)-10np\ln(20w))\right]=o(1).\label{eq:cprob}
				\end{align} 
		 
		Let $\mathcal{B}$ be the event that for every choice of $y\in [n]$ and $S\subset[n]\setminus y$ of size $m$, there exist $z_1 \in N(y) \setminus S$ and $z_2 \in[n] \setminus (N(y) \cup S \cup \{y\})$ that both do not have neighbours in $S$. Again, $\mathcal{B}$ implies that $\fun(G_n)\geq m$. However $(\neg\mathcal{B})\cap \mathcal{A}_{\mathrm{large}}  \subseteq \mathcal{C}$, and so $\mathbb{P}(\neg \mathcal{B}) =o(1)$ by \eqref{eq:cprob} and  Claim~\ref{clm:lower1}, as desired.   \end{proof}

			\subsection{Lower Bound for \texorpdfstring{$p\geq p^*$}{p at least p*}}\label{sec:lowerlargep} 
						
			\begin{lemma}\label{lem:lower-p-more-p-star}
			Let $G_n \sim G(n,p)$, where $p=\sqrt{\frac{w\ln n}{n}}\leq\frac{1}{2}$, $w=w(n)\geq 1$. Then  \whp
				\[
				\fun\left(G_n\right) \geq  \frac{\ln (ew)}{4p}. 
				\] 
			\end{lemma}
			\begin{proof}
				It is sufficient to prove that the following event $\mathcal{B}$ holds \whp: for every set $S$ of size 
				$$
				k:=\left\lceil\frac{\ln (ew)}{4p}\right\rceil>\frac{1}{2}(\ln n-\ln^2\ln n)
				$$ 
				and every vertex $y\notin S$, there exist two vertices $z_1,z_2\in[n]\setminus (S\cup\{y\})$ with no neighbours in $S$, where $z_1$ is adjacent to $y$, while $z_2$ is non-adjacent to $y$. Note that $\neg {\mathcal{B}}\subseteq \mathcal{E}_1 \cup \mathcal{E}_2$, where
				\begin{enumerate}
					\item [$\mathcal{E}_1$] is the event that there exist $S\subset[n]$ of size $k$ and  $y\notin S$ such that every neighbour of $y$ in $[n]\setminus (S\cup\{y\})$ has a neighbour in $S$;
					
					\item [$\mathcal{E}_2$] is the event that there exist $S\subset[n]$ of size $k$ and  $y\notin S$ such that every non-neighbour of $y$ in $[n]\setminus (S\cup\{y\})$ has a neighbour in $S$.
				\end{enumerate} 
				
				\paragraph{Bound on $\Pr{\mathcal{E}_1}$.}
				Let us fix $S\subset[n]$ of size $k$ and a vertex $y\notin S$. For every $z_1\in[n]\setminus (S\cup\{y\})$, the probability that $z_1$ is adjacent to  $y$ and $N(z_1)\cap S=\emptyset$ is exactly $p(1-p)^{k}$. Then
				$$
					\Pr{\nexists z_1\in[n]\setminus (S\cup\{y\})\,:\, N(z_1)\cap (S\cup\{y\})=\{y\}}=(1-p(1-p)^k)^{n-k-1},
				$$
				and by the union bound
				\begin{equation}
					\Pr{\mathcal{E}_1} \leq 
					\frac{n^{k}}{k!}\cdot n\cdot (1-p(1-p)^k)^{n-k-1} \leq
					\exp\left[k\ln n-(1/2-o(1))pn(1-p)^k\right].
					\label{eq:z_1}
				\end{equation}
				
				\paragraph{Bound on $\Pr{\mathcal{E}_2}$.}
				Let $\mathcal{R}$ be the event that every vertex has at least $n-np(1+o(1))\geq n(1/2-o(1))$ non-neighbours, and note that $\mathbb{P}(\neg\mathcal{R})=o(1)$ by Lemma~\ref{cl:degrees_common_degrees}. 	Fix $S'\in{[n]\choose k}$ and $y'\notin S'$. 
				Then, by the union bound and since $k!>n$, for sufficiently large $n$, $\Pr{\mathcal{E}_2}$ is bounded from above by
				\begin{align*}
					 & \Pr{\exists S\in{[n]\choose k}\,\exists y\notin S\,\forall z_2\notin(S\cup\{y\})\,:\, N(z_2)\cap (S\cup\{y\})\neq \emptyset } \\
					&\leq 
					{n\choose k}\cdot n\cdot\Pr{(n-|N[y']|\geq n(1/2-o(1)))\wedge(\forall z_2\in[n]\setminus (N[y']\cup S)\,:\, N(z_2)\cap S'\neq \emptyset)}+\mathbb{P}(\neg\mathcal{R})\\
					&\leq 
					n^k\cdot \Pr{\mathrm{Bin}(n(1/2-o(1)),(1-p)^k)=0}+o(1).
				\end{align*}
				Thus we have 
				\begin{equation}
					\Pr{\mathcal{E}_2} \leq n^k\left(1-(1-p)^k\right)^{n(1/2-o(1))}+o(1)
					\leq\exp\left[k\ln n- (1/2-o(1))n(1-p)^k\right]+o(1).
					\label{eq:z_2}
				\end{equation}
				
			\paragraph{Bound on $\Pr{\neg {\mathcal{B}}}$.}
			From~\eqref{eq:z_1}~and~\eqref{eq:z_2}, we get that
				\begin{equation}\label{eq:desiredbdd}
					\Pr{\neg {\mathcal{B}}} \leq \Pr{\mathcal{E}_1} + \Pr{\mathcal{E}_2}\leq 2\exp\left[k\ln n-(1/2-o(1))pn(1-p)^k\right]+o(1).\end{equation}
				We now show that the exponent on the right-hand side of \eqref{eq:desiredbdd} goes to $-\infty$ as $n$ goes to $+\infty$. 
				\begin{claim}\label{clm:lower2}Let $\rho(n):= k\ln n-(1/2-o(1))pn(1-p)^k$. Then, $\rho(n) =-\omega(1)$.
				\end{claim}
				\begin{pocd}{\ref{clm:lower2}}

					Set $\delta:=1/10$. First, let $ n^{-\delta}\leq p\leq 1/2$. Since $\varphi(p):=\ln(1-p)+2p\ln 2$ increases on $(0,1-1/(2\ln 2))$, decreases on $(1-1/(2\ln 2),1/2)$, $\varphi(0)=0$, and $\varphi(1/2)=0$, $\varphi(p)$ is non-negative on $[0,1/2]$ and thus $\frac{\ln(1-p)}{2p}\geq-\ln 2$ for $p\in[n^{-\delta},1/2]$. 
					Therefore, 
					\begin{align*}
						\rho(n)  \leq & 
						n^{\delta}\ln^2 n-(1/2-o(1))n^{1-\delta} \cdot \exp\left[\left(\frac{\ln (ew)}{4p}+1\right)\cdot\ln(1-p)\right] \\
						\leq  &
						n^{\delta}\ln^2 n-(1/2-o(1))n^{1-\delta} \cdot \exp\left[-\frac{\ln 2}{2}\cdot\ln\left(ew\right) + \ln(1-p)\right] \\
						\leq  &
						n^{\delta}\ln^2 n-n^{1-\delta} \cdot \exp\left[-\left(\frac{\ln 2}{2}+ o(1)\right) \ln n\right] \\
						\leq  &
						n^{2\delta}-n^{1-\delta-\ln 2/2-o(1)} =-\omega(1),
					\end{align*}
				for suitably large $n$.  
				
				Now, if $ \sqrt{\frac{\ln n}{n}}\leq p < n^{-\delta} $, then,  
				\begin{equation*}
					\rho(n)  \leq  (1+o(1))\frac{\ln(ew)}{4p}\ln n-\frac{1}{2}np\exp\left[-(1/4+o(1))\ln(ew)\right]. 					
				\end{equation*}
				Then recalling that $w= np^2/(\ln n)$ and noting that  $\ln x-(2/e)\cdot x^{5/7}$ decreases on $[e,\infty)$, we have
				\begin{align*}
					\rho(n) 	\leq & 
					\frac{ (1+o(1))\ln(ew)-2w\cdot (ew)^{-1/4-o(1)}}{4p/\ln n} \\
					\leq &
					\frac{ (1+o(1))\left(\ln(ew)-2w\cdot (ew)^{-2/7}\right)}{4p/\ln n}  \\
					\leq &  \frac{(1+o(1))\left(\ln e-(2/e)\cdot e^{5/7}\right)}{4p/\ln n} \\
					\leq &
					-\frac{(1+o(1))\left(n^{\delta} \cdot \ln n \right)}{8} 
					=-\omega(1),
				\end{align*}
				as claimed. 
				\end{pocd}
				Finally, by \eqref{eq:desiredbdd} and \cref{clm:lower2}, we have $\Pr{\neg {\mathcal{B}}} =o(1)$, which completes the proof. 
			\end{proof}

			\section{Domination in Random Bipartite Graphs}
			\label{sc:bip_dom}

			Let $G=(V,E)$ be a graph and $S \subseteq V$. We say that a set $D \subseteq V \setminus S$ is a \emph{dominating set} for $S$ if every vertex in $S$ has a neighbour in $D$. Recall that for two independent sets $A,B$ and $p\in[0,1]$, the random bipartite graph $G(A,B,p)$ denotes the distribution over bipartite graphs on $A,B$ where each edge $ab$ with $a\in A$ and $b\in B$, is present independently with probability $p$.
			
			The following result (\cref{lem:small_t_dom_exist}) bounds dominating sets in the binomial random graph, and, as mentioned in the introduction, this bound is tight. \cref{lem:small_t_dom_exist} is a key ingredient in the proof of the upper bound on functionality, and as it is written primarily with this task in mind, the statement is a little opaque. After proving \cref{lem:small_t_dom_exist} we will show in \cref{sec:simpleDomLemma} how \cref{lem:small_t_dom_simple}, the simplified version of \cref{lem:small_t_dom_exist}, follows from \cref{lem:small_t_dom_exist}.

			\DomLemma*
			\begin{proof} 
				We will now sketch the proof, which uses a greedy algorithm to find a subset of vertices in $B$ that dominates most of $A$. To begin, we assign an arbitrary ordering $\sigma$ to the vertices of $B$. The algorithm then proceeds in rounds,  where in round $i\geq 1$, for ``thresholds'' $k_i$ and $\ell_i$ to be defined, we expose the edges from vertices in $B$ vertex-by-vertex, following $\sigma$, and we add $v\in B$ to the dominating set $D$ if it has more than $k_i$ neighbours in $A$ that are not yet dominated. The round ends when there are less than $\ell_i$ vertices of $A$ that are not yet dominated. The thresholds $k_i$ and $\ell_i$ are both non-increasing in $i$, but $k_i$ decays at a slower rate than $\ell_i$. This algorithm terminates when there are at most $1/(ph)$ vertices of $A$ which are not dominated, or if we have exposed all edges incident with vertices in $B$. Assuming we do not exhaust $B$, the set $D$ will contain the desired number of vertices by definition of $k_i$ and $\ell_i$. Thus, if the algorithm terminates before we have ``run out'' of vertices in $B$, then we have achieved our goal. Our task is to bound the probability of this event. 
				
				Notice that by the conditions $\ln^2 n\leq a$ and $aph \leq \alpha \ln n  $ we have	
				\begin{equation}\label{eq:condpw} 
					ph \leq \frac{\alpha \ln n }{\ln^2 n}  = \frac{\alpha}{\ln n}.  
				\end{equation}
				Further, observe that if $a \leq 1/ph$, then we can set $A'=A$ and the theorem holds trivially. Thus, we can assume that		 
				\begin{equation}\label{eq:aphbig1}
					aph >1. 
				\end{equation}
				We now set, with foresight, the constant  
				
				\begin{equation}\label{eq:cchoice} c:= \max\left\{ \frac{20\alpha}{\delta} , \;1\right\},\end{equation} and define the ``thresholds'' $k_i,\ell_i$ in  the following way: $\ell_0:=a$, and, for $i\geq 1$, 
				\begin{equation}
					\label{eq:kldef} 
					k_i : = \left\lceil \frac{aph}{i  c^2 \ln(2h)}  \right\rceil \qquad \text{and} \qquad \ell_i:= \left\lceil a \cdot e^{-i\cdot c}\right\rceil.  
				\end{equation}
				Consider the setting where, given our ordering $\sigma$ of $B$, we extend $\sigma$ to an ordering of a countably infinite set $U$ --- thus $B$ is set of  the first $b$ elements of $U$ in this ordering. We are then considering running the algorithm on a random bipartite graph $G(A,U,p)$. Before starting the algorithm, colour each vertex of $U$ green. In the $i$-th round of the algorithm, where $i\geq 1$, we initialise $X_i=\emptyset$. We then process vertices one-by-one by choosing the first green vertex from $U$ in the ordering $\sigma$, say $v$, and first colour it red. We then expose $N(v)$ and add $v$ to $X_i$ if $|N(v)\setminus N(D_i)| \geq k_i$, where $D_i= \bigcup_{1\leq j \leq i}X_j$.  The $i$-th round ends when $|A\setminus N(D_i)|< \ell_i$, and we let $\tau_i\geq 0$ denote the number of vertices in $U$ which changed colour from green to red during the $i$-th round. The algorithm ends after some (random) number $T$ of rounds such that $|A\setminus N(D_T)| \leq 1/(ph)$, and we let $D:=D_T$.
				
				As we shall see shortly (\cref{clm:k-i-l-i}), $k_i\leq \ell_i$ holds for large $n$ and all $i \in [T]$ almost surely, and thus, since $p>0$,
				this algorithm terminates almost surely. We say that the algorithm is successful if the number of vertices of $U$ which are turned red is at most $b$, that is the event $\{\sum_{i=1}^T\tau_i \leq b\}$ holds. In this case we find the desired set $D\subseteq B$. 
				
				We now prove three bounds, one for $T$, and two which relate $k_i$ to $\ell_i$. To begin, observe that  \begin{equation}\label{eq:boundsl_t}
					\frac{1}{2e^{c}ph}=
					a\cdot e^{-((1/c) \cdot \ln(2aph) + 1)\cdot c}
					\leq a\cdot e^{-\lceil (1/c) \cdot \ln(2aph)\rceil\cdot c}\leq 
					\Big\lceil a\cdot e^{-\lceil (1/c) \cdot \ln(2aph)\rceil\cdot c} \Big\rceil
					\leq
					\Big\lceil \frac{1}{2ph} \Big\rceil
					\leq \frac{1}{ph}, 
				\end{equation} since $1/(ph) = \Omega(\ln n) $ by \eqref{eq:condpw}, $aph>1$ by \eqref{eq:aphbig1}, and $n$ is suitably large. Thus, from the definition of $T$ and the definition of $\ell_i$ given in~\eqref{eq:kldef}, it follows that    \begin{equation}\label{eq:boundonT} 
					T\leq T':=  \lceil (1/c) \cdot \ln(2aph)\rceil\quad\text{almost surely}. 
				\end{equation} 
				Therefore, almost surely, for all $1\leq i \leq T$, we have  $\ell_i =\omega(1)$ and so from \eqref{eq:kldef} we obtain
				\begin{equation}\label{eq:klrel}
					p\cdot \ell_i\leq p\cdot\ell_{T'}  
					\leq  (1+o(1))\cdot k_i \cdot \frac{i  c^2\ln(2h) }{h}\cdot    e^{-i\cdot c}. 
				\end{equation}
				
				\begin{claim}\label{clm:k-i-l-i}
					For any $n$ large enough and all $i \in [T]$, we have $k_i\leq \ell_i$ almost surely.
				\end{claim}
				\begin{pocd}{\ref{clm:k-i-l-i}}
					Note that  since $c\geq 1$ by \eqref{eq:cchoice},  it follows that $i\cdot e^{-i\cdot c}$ is decreasing in $i\geq 1$. Thus by \eqref{eq:boundonT}, almost surely, for all $1\leq i\leq T$
					\[ \frac{aph}{i  c^2 \ln(2h)}  \cdot \frac{1}{ a \cdot e^{-i\cdot c}}  \leq \frac{aph}{T  c^2 \ln(2h)}  \cdot \frac{1}{ a \cdot e^{-T\cdot c}} \leq \frac{aph}{\lceil (1/c) \cdot \ln(2aph)\rceil\cdot   c^2 }  \cdot \frac{1}{ a\cdot e^{-\lceil (1/c) \cdot \ln(2aph)\rceil\cdot c} }. \]
					Plugging in the bound from \eqref{eq:boundsl_t}, $aph \leq \alpha \ln n$, and $a\geq \ln^2 n$, for large $n$ we have
					\begin{equation}\label{eq:kleql} 		
						\frac{aph}{i  c^2 \ln(2h)}  \cdot \frac{1}{ a \cdot e^{-i\cdot c}}			\leq \frac{aph}{ \ln(2aph)  }  \cdot 2e^{c}ph  
						= \frac{2e^{c}(aph)^2}{  \ln(2aph)\cdot a   }   
						\leq \frac{2e^{c}(\alpha \ln n )^2}{  \ln(2\alpha \ln n ) \cdot \ln^2 n   }   \leq 1,
					\end{equation}
					where the penultimate inequality uses the fact that $x^2/\ln(2x)$ is increasing on the interval $x > 1$.
					Now, since $x\leq y$ implies $\lceil x\rceil \leq \lceil y\rceil$, it follows from \eqref{eq:kleql} that  $k_i \leq \ell_i$ for large $n$ and all $i\in [T]$ almost surely.
				\end{pocd}
				
				Observe that by the definition of the rounds we only have an upper bound on how many vertices of $A$ are not dominated at the start of each round, also it may happen that some round $i$ is actually skipped (i.e.~$\tau_i=1$) if we find a vertex in a previous round that dominates many vertices of $A$. 
				Let $\mathcal{E}_{i,j}$ be the event that the $j$-th vertex uncovered during the $i$-th round has at least $k_i\leq \ell_i$ neighbours among the un-dominated vertices in $A$, and let $\mathcal{B}_{i,j}$ be the event that, at the time when the algorithm processes the $j$-th vertex in the $i$-th round we still have that the un-dominated part of $A$ has cardinality at least $\ell_i$. Then 
				\begin{equation} \label{eq:lowerEi}
					\Pr{\mathcal{E}_{i,j}\mid\mathcal{B}_{i,j}} \geq \Pr{\bin{\ell_i}{p} = k_i} = \binom{\ell_i}{k_i}p^{k_i}(1-p)^{\ell_i-k_i}\geq \left(\frac{ p\cdot \ell_i }{k_i}  \right)^{k_i}e^{-2\cdot p\cdot \ell_i},  \end{equation}
				where the last inequality holds since $e^{-2p} \leq 1-p$ holds for all $p \in [0,.79]$, and $p\leq ph\leq \alpha/\ln n$ by \eqref{eq:condpw}, and $n$ is large by assumption. We now bound the right-hand side of~\eqref{eq:lowerEi} from below. Let $p_i:=\Pr{\bin{\ell_i}{p} = k_i}$.
				
				\begin{claim}\label{clm:keydominating} For any $1\leq i\leq T'$, we have  $p_i\geq  e^{-2c-1} \cdot n^{-\delta}$. 
				\end{claim}
				\begin{pocd}{\ref{clm:keydominating}} We have two cases for different values of $h$.

					The first case is $\ln(2h) \geq \frac{aph}{ic^2}$. Observe that in this case $k_i=1$ holds by \eqref{eq:kldef}. Thus by \cref{eq:klrel}, since $x/e^x$ is maximized by $x=1$ and $\ln(2x)/x$ is maximized by $x=e/2$, for large $n$ we have
					\[	
					p\cdot \ell_i   
					\leq  (1+o(1))\cdot k_i \cdot \frac{i  c^2}{e^{i\cdot c}} \cdot \frac{\ln(2h) }{h} 
					\leq (1+o(1)) \cdot 1\cdot \frac{c}{e} \cdot \frac{\ln(e) }{e/2} 
					\leq \frac{c}{2}.
					\] 
					Observe also that $  p \cdot \ell_i \geq  p \cdot \ell_{T'}\geq  1/(2e^{c}h)$ by \eqref{eq:boundsl_t}. Thus, by \eqref{eq:lowerEi} and $h\leq n^\delta$, and recalling that $k_i=1$, we have 
					\[
					p_i \geq   (p\cdot \ell_i) \cdot  e^{-2\cdot p\cdot \ell_i} \geq \frac{1}{2e^{c}h} \cdot e^{-c} \geq e^{-2c-1} \cdot n^{-\delta}.  
					\] 
					
					The second case is $\ln 2\leq \ln(2h) < \frac{aph}{ic^2}$. This condition  implies that $\frac{aph}{ic^2\ln(2h)}>1$, and thus, $k_i\geq 2$. Therefore, 		 
					\begin{equation*}	
						p\cdot \ell_i 
						\geq  
						p \cdot a \cdot e^{-i\cdot c}  
						= 
						\frac{aph}{i  c^2 \ln(2h)}\cdot \frac{i  c^2\ln(2h)}{h}  \cdot e^{-i\cdot c}
						>  
						\frac{k_i}{2} \cdot 	\frac{i  c^2\ln(2h) }{h}\cdot    e^{-i\cdot c}. 
					\end{equation*}	
					However, for large $n$, we have $
					p\cdot \ell_i  
					\leq  2 k_i \cdot \frac{i  c^2\ln(2h) }{he^{i\cdot c}}$  by \eqref{eq:klrel}. Using these bounds in  \eqref{eq:lowerEi} gives
					\begin{align}
						p_i &\geq  \left(\frac{1}{2}\cdot\frac{     i c^2\ln(2h) }{h} \cdot e^{-i\cdot c } \right)^{k_i}\cdot \exp\left[-2\cdot  2 k_i \cdot\frac{ i c^2\ln(2h) }{he^{i\cdot c}}  \right] \notag \\
						&=\exp\left[ -\left(   i  c\left(1+\frac{4c\ln(2h)}{he^{i\cdot c }}\right)  - \ln\left(\frac{ i c^2\ln(2h) }{2h } \right) \right) \cdot k_i \right].\label{eq:boundPei}
					\end{align}
					
					From the conditions $i,h,c\geq 1$ we have $\ln(2h) < h$ and $ce^{-i\cdot c} < 1/2$, which give the bounds
					\[ \frac{4c\ln(2h)}{he^{i\cdot c }} < 2\qquad \text{and} \qquad -\ln\left(\frac{ i c^2\ln(2h) }{2h }\right) = \ln\left(\frac{2h }{ i c^2\ln(2h) }\right)
					\leq  \ln\left(\frac{2h}{\ln 2}   \right) < \ln\left(3h   \right)  .\]
					Inserting these bounds into \eqref{eq:boundPei}  yields

					\begin{equation*} 
						p_i>  \exp\left[ -2  \cdot\left(   3i  c   + \ln(3h) \right) \cdot \frac{aph}{i  c^2 \ln(2h)} \right]	=    \exp\left[ - \frac{6aph}{c\ln(2h)}  -  \frac{2aph\ln (3h)}{ic^2 \ln (2h)}   \right].
					\end{equation*}
					Using   $i,h,c\geq 1$, $1/\ln (2h) \leq 3/2$, $\ln(3h)/\ln(2h)\leq \ln(3)/\ln(2)< 2$, and $aph \leq \alpha \ln n$ gives 
					\begin{equation*} 
						p_i>      \exp\left[ - \frac{9aph}{c}   -  \frac{4aph}{c }   \right]		 	>  n^{- 13\alpha/c } > n^{- \delta } ,  
					\end{equation*}where the last inequality holds by our choice of $c$ from \eqref{eq:cchoice}. 
				\end{pocd}

				We now show that the ``success'' event $\{\sum_{i=0}^{T}\tau_i \leq  b\}$ holds with the desired probability. First, observe that in round $i\geq 1$, there are at most $\ell_{i-1}-1 $ vertices of $A$ which are not dominated at the start of round $i$,  when each vertex is added to $X_i$ it covers at least $k_i$ undominated vertices in~$A$, and before the last vertex of round $i$ is added to $X_i$ there are at least $ \ell_i$ undominated vertices in~$A$. Thus, we have \[\ell_{i-1} -1 - (|X_i|-1)\cdot k_i \geq \ell_{i}.\] This, together with \eqref{eq:klrel}, implies that for $1\leq i\leq T$, almost surely, we have  
				\[|X_i| \leq \frac{\ell_{i-1} -1 - \ell_{i}}{k_i} + 1 \leq  \frac{\ell_{i-1}}{k_i} +1 \leq  (1+o(1))e^c\cdot \frac{\ell_{i}}{k_i} +1  \leq (1+o(1))\cdot  \frac{ c^2e^c\ln(2h) }{ph}\cdot   i  e^{-i\cdot c} + 1.\] 
				Thus, as $c\geq 1$ by \eqref{eq:cchoice}, $\sum_{i=1}^\infty ie^{-i} = \frac{e}{(e-1)^2}<0.93$, $1/ph =\Omega(\ln n)$ by \eqref{eq:condpw}, $T =O(\ln \ln n )$ by \eqref{eq:boundonT},  and $n$ is large, we have 
				\begin{equation}\label{eq:sizeofXi} 
					\sum_{i \in [T]}|X_i| \leq T +  (1+o(1))\cdot  \frac{ c^2e^c\ln(2h) }{ph}\cdot \sum_{i \in [T]} i e^{-i\cdot c} \leq   \frac{ c^2e^c\ln(2h) }{ph}\quad\text{ almost surely}.     
				\end{equation}

				Due to~\eqref{eq:lowerEi}, there exist independent geometric random variables $\eta_j^i\sim \geo{p_i}$,  $j\in [|X_i|]$, such that, almost surely, for each $i\in[T]$,  $\tau_i\leq\sum_{j=1}^{|X_i|}\eta_j^i$.   
				Let 
				\begin{equation}\label{eq:dfp*} 
					p_* := e^{-2c-1} \cdot n^{-\delta}  \quad\text{ and }\quad N:=   \left\lfloor\frac{c^2e^c\ln(2h) }{ph}\right\rfloor.
				\end{equation} 		
				We know that, almost surely, $\min_{ i \in [T]}p_i \geq p_*$ by \cref{clm:keydominating}, and  $\sum_{i \in [T]}|X_i| \leq  N$ by \eqref{eq:sizeofXi}. Hence there exist independent $\zeta_k\sim \geo{p_*}$ for $k\in [N]$ such that
				\[
				\sum_{i \in [T]}\tau_i \leq  \sum_{i \in [T] } \sum_{j=1}^{|X_i|} \eta_{j}^i  \leq  \sum_{k=1}^{N} \zeta_k\quad\text{ almost surely}.
				\]
				Thus, if $Z:=\sum_{k=1}^N  \zeta_k $, then it suffices for our goals to bound $Z$ from above. To begin, observe that by \eqref{eq:dfp*}, 
				\begin{equation}\label{eq:mup*}\mu:=\Ex{Z}= \frac{1}{p_*}\cdot N = e^{2c+1}   n^{\delta}\cdot  \left\lfloor\frac{c^2e^c\ln(2h) }{ph}\right\rfloor.
				\end{equation} 
				By assumption we have $bp>n^{\beta}$, $h \leq n^{\delta}$, and $\delta<\beta$ are constants, so by \eqref{eq:mup*} we have $b/\mu =\omega(1)$. Thus, for large $n$, by applying \cref{lem:jansontail} with $\lambda = b/\mu$, and $p_*=e^{-2c-1} \cdot n^{-\delta}$, we have  
				\begin{align*}
					\Pr{\sum_{i \in [T]}\tau_i > b  } 
					& \leq 
					\Pr{Z > b}   
					\\ & \leq 
					\exp\left[-  p_*\cdot \mu\cdot (\lambda - 1 -\ln \lambda)\right] 
					\\ & \leq 
					\exp\left[-(1-o(1))\cdot p_*\cdot b\right] 
					\\ & \leq 
					\exp\left[-e^{-2c-2}\cdot n^{-\delta}\cdot b\right]. 
				\end{align*}
				The result follows by taking $C:=\max\{c^2e^c,\, e^{2c+2}\}$ in the statement, as conditional on $\sum_{i \in [T]}\tau_i\leq b$  the desired set $D:=D_T$ has size at most $ \sum_{i \in [T]}|X_i|   <  \frac{C\ln(eh) }{ph}$ almost surely by \eqref{eq:sizeofXi}.  \end{proof}

			\subsection{Simplified Version of \cref{lem:small_t_dom_exist}}\label{sec:simpleDomLemma}
			
			In this section, we show how to derive from \cref{lem:small_t_dom_exist} its simplified version stated in the introduction, which we restate below for convenience.

			\simpleDomLemma*
			\begin{proof}
				Fix $h=1$ in the \cref{lem:small_t_dom_exist} and thus $h\leq n^\delta$ for any $\delta\geq 0$. For $\alpha, \beta >0$ in the statement of \cref{lem:small_t_dom_simple} and  any $\delta >0$, let $\mathcal{D}_\delta$ denote the event in the statement of \cref{lem:small_t_dom_exist}, and let $\mathcal{C}:=\{\forall a\in A,\; N(a)\neq \emptyset\} $. Observe that $\mathcal{C}$ and $\mathcal{D}_\delta$ are both monotone increasing properties. Thus, as a consequence of the FKG inequality \cite[Theorem 6.3.3]{AlonSpencer}, we have $\Pr{\mathcal{D}_\delta\mid \mathcal{C}} \geq \Pr{\mathcal{D}_\delta} $. Conditional on the event $\mathcal{C}$, the set $A'$ described in the event $\mathcal{D}_\delta$ can be dominated using at most $1/p$ additional vertices from $B$, giving a dominating set $D'\subseteq B$ of $A$ with size at most $C'/p$ where $C':= 1+C(\alpha, \delta)\ln 2$. Now, since $f(n)\rightarrow\infty$, it follows that for any $\varepsilon >0$, there  exists constants $N,\delta_N>0$ such that for $n\geq N$ we have $n^{-\delta_N} b/C(\alpha,\delta_N) \geq b^{1-\varepsilon} $ and $C'=(1+C(\alpha,\delta_N) \ln 2)\leq f(n)$. The result then follows by applying \cref{lem:small_t_dom_exist} over such a sequence of $ \delta_N$'s. 
			\end{proof}

			We conclude this section with two propositions showing the optimality of \cref{lem:small_t_dom_simple}. Namely, \whp~a minimum dominating set of $A$ in $B$ has size $\Omega(1/p)$ when $ap=\Theta(\ln n)$ (\cref{prop:dom-lower-bound}), and the probability of existence of a dominating set of size $O(1/p)$ is at most $1-\exp\left[-b^{1+o(1)}\right]$ when $ap=\Omega(1)$ (\cref{prop:dom-upper-bound}).
			
			\begin{proposition}\label{prop:dom-lower-bound}
				Let $a,b$ and $p$ satisfy the conditions of \cref{lem:small_t_dom_simple}, and additionally $ap \geq c\ln n$ for some constant $c>0$. Then, for large $n$,  	\[	\Pr{ \neg\exists \text{ a dominating set $D\subset B$ of }A\text{ with } |D|\leq \frac{c}{2p}   } 
				\geq  1- \exp\left[- a/3 \right]. 
				\]
			\end{proposition}
			\begin{proof}
				We can assume w.l.o.g.~that $c\leq 2/e$. Let $X$ be the number of dominating sets $D\subseteq B$ of $A$ which satisfy $|D|\leq d:=  \frac{c}{2p} \leq  \frac{1}{ep} $. Then, 
				\[\Ex{X} \leq d \cdot \binom{b}{d}\cdot (1-(1-p)^d)^a \leq b\cdot b^d\cdot (dp)^a \leq n \cdot e^{d\ln n}\cdot e^{-a} \overset{(i)}{\leq} n \cdot  e^{-a/2} \overset{(ii)}{\leq} e^{-a/3}, \] where $(i)$ holds as $d\ln n \leq \frac{c\ln n}{2p} \leq a/2 $ and $(ii)$ holds for large $n$ since $a \geq \ln^2 n$ by the hypothesis from \cref{lem:small_t_dom_simple}.  The result then follows by Markov's inequality.
			\end{proof}

			\begin{proposition}\label{prop:dom-upper-bound}
				For any constant $C\geq 1$ and $a,b$ and $p$ satisfying the conditions of \cref{lem:small_t_dom_simple}, and additionally $ap > C$ 
				\[	\Pr{ \exists \text{ a dominating set $D\subset B$ of }A\text{ with }\left. |D|\leq \frac{C}{p} \, \right| \, \forall a\in A,\; N(a)\neq \emptyset } 
				\leq  1- \exp\left[- 3abp \right]. 
				\]
			\end{proposition}
			\begin{proof}
				Choose any $\ln^2 n \leq a \leq b \leq n $ and $p\in [0,1]$ which satisfy $bp\geq n^{\beta}$.	Let \[\mathcal{E}:=\{\exists \text{ a dominating set $D\subset B$ of }A\text{ with } |D|\leq C/p\},\quad \text{and}\quad \mathcal{C}:=\{\forall a\in A,\; N(a)\neq \emptyset\}, \]be the two events in the statement, and further let $\mathcal{M}:=\{E(G)\text{ is a matching of }A\text{ in }B\}$. Observe that since $a>C/p$ we have $\mathcal{M} \subseteq \neg \mathcal{E} $. Consequently, as $\Pr{\mathcal{C}}>0$, we have 
				\[\Pr{\mathcal{E} \mid \mathcal{C}} = 1- \Pr{\neg \mathcal{E} \mid \mathcal{C}} \leq 1 - \Pr{ \mathcal{M} \mid \mathcal{C}} = 1 - \frac{\Pr{ \mathcal{M}} }{ \Pr{ \mathcal{C}}} \leq 1- \Pr{ \mathcal{M}}.  \] 
				Now since there are $\binom{b}{a}\cdot a! $ ways to choose a matching of $A$ in $B$ we have 
\begin{align*}
\Pr{ \mathcal{M}} & = \binom{b}{a}a!\cdot p^a(1-p)^{(b-1)a}\\
 &\geq (b-a)^a p^a (1-p)^{ab}\\
 & \geq b^a e^{-2a^2/b} \cdot p^a e^{-2abp} =\exp[a\ln(bp)  - 2a^2/b - 2abp ].
\end{align*}
Recalling that $n^\beta \leq bp \leq n$ and observing $2a^2/b \leq abp$ gives $\Pr{ \mathcal{M}} \geq \exp[- 3abp ]$.   
			\end{proof}

\section{Conclusions}

In this paper we determined the functionality of $G(n,p)$ up to a constant factor for all $p\in[0,1]$. In particular, \whp the functionality of $G(n,p)$ is at most $\Theta(\sqrt{\frac{n}{\ln n}})$ for any $p\in [0,1]$ (cf. \cref{fig:fun-plot}).
This is lower than the largest known functionality of an $n$-vertex graph, which is $\Omega(\sqrt{n})$ and is achieved on the incidence graph of a finite projective plane \cite{DFOPSA23}. On the other hand, the best known general upper bound on the functionality of an $n$-vertex graph is $O(\sqrt{n\ln n})$ \cite{DFOPSA23}.
We believe that this general upper bound can be improved, and suspect that it should be possible to do this all the way to $O(\sqrt{n})$.

\begin{problem}
	Does there exist a constant $C$ such that the
	functionality of every $n$-vertex graph is at most $C\sqrt{n}$?
\end{problem}

\noindent
In order to give some evidence supporting our intuition, we first recall that both upper and lower bounds from \cite{DFOPSA23} utilize the observation
that in a graph $G=(V,E)$, a vertex $v$ is a function of a set $S$ if and only if $S$ intersects every set $(N(u_1) \triangle N(u_2)) \cup \{u_1,u_2\}$, where $u_1$ is a neighbour and $u_2$ is a non-neighbour of $v$, and $u_2 \neq v$. In other words, $v$ is a function of $S$ if and only if $S$ is a transversal in the hypergraph $\mathcal{H}(v) := (V, \{ (N(u_1) \triangle N(u_2)) \cup \{u_1,u_2\}~:~ u_1 \in N(v), u_2 \not\in N[v] \})$.

Using this observation, the upper bound of $O(\sqrt{n\ln n})$ is obtained by showing that if, for every two vertices $u,v$, the symmetric difference $N(u)\triangle N(v)$ has size at least $\Omega(\sqrt{n\ln n})$, then a uniformly random set of size $\Theta(\sqrt{n\ln n})$  intersects each $N(u)\triangle N(v)$ with a positive probability. The lower bound of $\Omega(\sqrt{n})$ is proved by showing that in the incidence graph $P$ of a finite projective plane with $n$ points and $n$ lines, for every vertex $v$, no set of size less than $\sqrt{n/2}$ is a transversal in $\mathcal{H}(v)$.
 
 A possible way to improve the lower-bound construction is to add or delete some random edges in $P$. This is supported by the following two observations. First, while the functionality of $K_n$ is 0, after removing every edge from $K_n$ with probability $1/2$ the functionality of the resulting graph becomes $\fun(G(n,p))=\Theta(\log n)$ \whp Second, it is known that the minimum transversal of a hypergraph with suboptimal transversal number typically increases after deleting half of vertices from every hyperedge uniformly at random (see,~\cite[Theorem 4]{HY-arxiv}). However, this strategy does not seem to work for $P$. Indeed, the maximum degree (and thus functionality) of $P$ is $O(\sqrt{n})$ and the removal of edges can only lower the maximum degree. On the other hand, we might add at least $\Theta(n\sqrt{n\ln n})$ random edges to $P$ instead, but it makes the influence of deterministic edges insignificant. In particular, it can be shown that any distinguishing set $S$ in the original graph of size $\Theta(\sqrt{n})$, \whp is also a distinguishing set in the randomly perturbed graph. 

Finally, we note that there exist $k$-uniform hypergraphs on $n$ vertices with the transversal number $\Omega(n\log k/k)$ (see \cite{Alon,TY,HY-arxiv}), which would be promising constructions for attempting to lower bound  functionality with $\Omega(\sqrt{n\ln n})$ (by taking $k=\Theta(\sqrt{n\ln n})$). Unfortunately, it seems unlikely that for such a hypergraph there exists a corresponding graph with the symmetric differences coinciding with the hypergraph's hyperedges.

			\bibliographystyle{alpha}

			\bibliography{bibliography}

		\end{document}